\documentclass[10pt]{article}
\usepackage{graphicx}
\usepackage{epsfig}
\usepackage{amssymb}
\usepackage{amsmath}
\usepackage{amsfonts}
\usepackage{dsfont}
\newcommand{\RR}{\mathbb{R}}
\newcommand{\NN}{\mathbb{N}}

\def\div{\operatorname{div}}

\newtheorem{theorem}{Theorem}[section]

\newtheorem{corollary}[theorem]{Corollary}
\newtheorem{proposition}[theorem]{Proposition}
\newtheorem{e-definition}[theorem]{Definition}
\newtheorem{remark}{\bf Remark\/}

\renewcommand{\div}{\operatorname{div}}

\begin{document}

%\begin{frontmatter}

% Titre, auteurs et adresses

% utiliser la commande \thanksref dans \title, \author ou \address
%     pour les notes en bas de page ;

% utiliser la commande \ead pour l'adresse e-mail de chaque auteur
%    (aprËs la commande \auteur) ;

% \author{Name\thanksref{label2}}
%\selectlanguage{english}
\title{On correctors for linear elliptic homogenization in the presence of local defects: the case of advection-diffusion}
\author{X. Blanc$^1$, C. Le Bris$^2$ \& P.-L. Lions$^3$\\
{\footnotesize $^1$ Universit\'e Paris Diderot, Laboratoire Jacques-Louis Lions,}\\
{\footnotesize B\^atiment Sophie Germain, 5, rue Thomas Mann}\\
{\footnotesize 75205 Paris Cedex 13, FRANCE,}\\
 {\footnotesize{\tt blanc@ann.jussieu.fr}}\\
{\footnotesize $^2$ Ecole des Ponts and INRIA,}\\
{\footnotesize 6 \& 8, avenue Blaise Pascal, 77455 Marne-La-Vall\'ee
  Cedex 2, FRANCE}\\
{\footnotesize{\tt lebris@cermics.enpc.fr}}\\
{\footnotesize $^3$ Coll\`ege de France, 11, place Marcelin Berthelot,}\\
{\footnotesize 75231 Paris Cedex 05,  and}\\
{\footnotesize CEREMADE, Universit\'e Paris Dauphine, Place de Lattre de
  Tassigny,}\\
{\footnotesize 75775 Paris Cedex 16, FRANCE}\\
{\footnotesize \tt lions@ceremade.dauphine.fr}
}
\maketitle

%
%\vspace{-2.6cm}
%%\selectlanguage{francais}
%
%\title{Correcteurs en homog\'en\'eisation elliptique lin\'eaire en pr\'essence de d\'efauts locaux: le cas de l'advection-diffusion}
%
%\author[xavier]{Xavier Blanc}
%\ead{blanc@ann.jussieu.fr}
%\author[claude]{Claude Le Bris}
%\ead{lebris@cermics.enpc.fr}
%\author[pll]{Pierre-Louis Lions}
%\ead{lions@ceremade.dauphine.fr}
%
%\address[xavier]{Universit\'e Paris Diderot, Laboratoire Jacques-Louis Lions, B\^atiment Sophie Germain,\\ 5, rue Thomas Mann, 75205 Paris Cedex 13,}
%\address[claude]{Ecole des Ponts and Inria,  77455 Champs~sur~Marne Cedex\thanksref{claudethanks}}
%\thanks[claudethanks]{The work of the second author is partly supported by  ONR under Grant N00014-15-1-2777 and by EOARD.}%under Grant ?? 
%\address[pll]{Coll\`ege de France, 11, place Marcelin Berthelot, 75231
%  Paris Cedex 05, and\\ CEREMADE, Universit\'e Paris Dauphine,
%place du Mar\'echal de Lattre de Tassigny, 75775 Paris Cedex 16, FRANCE}
%% 
%

%\selectlanguage{english}

\begin{abstract}
%
%\noindent
%\vskip 0.5\baselineskip
%% Texte de l'abstract en anglais
\noindent We follow-up on our  works devoted to homogenization theory for  linear second-order  elliptic equations with coefficients that are perturbations of periodic coefficients. We have first considered  equations in divergence form in~\cite{milan,cras-defauts,cpde-defauts}. We have next shown, in our recent work~\cite{BLL-2017-1}, using a slightly different strategy of proof than in our earlier works, that  we may also address the equation~$-a_{ij}\partial_{ij}u=f$. The present work is devoted to advection-diffusion equations:~$-a_{ij}\partial_{ij}u+b_{j}\partial_ju=f$. We prove, under suitable assumptions on the coefficients $a_{ij}$, $b_{j}$, $1\leq i,j\leq d$ (typically that they are the sum of a periodic function and some perturbation in~$L^p$, for suitable $p<+\infty$), that the equation admits a (unique) invariant measure and that this measure may be used to transform the problem into a problem in divergence form, amenable to the techniques we have previously developed for the latter case.  
% \medskip
% %\selectlanguage{francais}
% \noindent{\bf R\'esum\'e}
% \vskip 0.5\baselineskip
% % Texte du r\'esum\'e en français
% \noindent Nous continuous ici notre s\'erie de travaux sur l'homog\'en\'eisation d'\'equations lin\'eaires elliptiques du second ordre dont les coefficients sont des perturbations de coefficients p\'eriodiques. Nous avons consid\'er\'e les \'equations sous forme conservative dans~\cite{milan,cras-defauts,cpde-defauts}, et, grâce \`a une strat\'egie de preuve l\'eg\`erement diff\'erente, l'\'equation~$-a_{ij}\partial_{ij}u=f$ dans notre travail plus r\'ecent~\cite{BLL-2017-1}. Ici, nous \'etudions les \'equations d'advection-diffusion~:~$-\partial_i(a_{ij}\partial_ju)+b_{j}\partial_ju=f$. Nous montrons, sous des hypoth\`eses ad\'equates sur les coefficients~$a_{ij}$, $b_{j}$, $1\leq i,j\leq d$ (typiquement qu'ils sont somme d'une fonction p\'eriodique et d'une fonction dans~$L^p$, pour $p<+\infty$ convenable), que l'\'equation admet une (unique) mesure invariante qui peut \^etre utilis\'ee pour transformer l'\'equation en une \'equation sous forme divergence, \'equation qu'on peut alors traiter par les techniques  d\'evelopp\'ees pr\'ec\'edemment pour ce cas. 
\end{abstract}

\newpage
%%%%%%%%%%%%%%%%%%%
\section{Introduction}
\label{sec:intro}

We study homogenization theory for the advection-diffusion equation
\begin{equation}
  \label{eq:equation-adv-diff-epsilon}
    -\, a_{ij}(x/\varepsilon)\,\partial_{ij} u^\varepsilon +\varepsilon^{-1}\,b_j(x/\varepsilon)\,\partial_{j} u^\varepsilon =f,
\end{equation} 
when the coefficients $a$ and $b$ in~\eqref{eq:equation-adv-diff-epsilon}  are \emph{perturbations}, formally
vanishing at infinity, of  periodic coefficients. Equation~\eqref{eq:equation-adv-diff-epsilon} is supplied with homogeneous Dirichlet boundary conditions and posed on a bounded regular domain $\Omega\subset\RR^d$,  with a right-hand-side term~$f\in L^2(\Omega)$. We assume that the coefficients $a$ and $b$ satisfy
\begin{equation}
  \label{eq:aper+tildea}
  a=a^{per}+\tilde a, \quad b=b^{per}+\tilde b
\end{equation}
where $a^{per}$, $b^{per}$ describe a \emph{periodic} unperturbed background, and $\tilde a$, $\tilde b$ the perturbation, with 
\begin{equation}
\label{eq:hyp1}
\left\{
\begin{array}{l}
a^{per}(x)+\tilde a(x) \quad\hbox{\rm and}\quad a^{per}(x) \quad \hbox{\rm are both uniformly elliptic, in}\,x\in\RR^d, 
\\
a^{per}\in \left(L^\infty(\RR^d)\right)^{d\times d},\quad b^{per}\in \left(L^\infty(\RR^d)\right)^d,
\\
  \tilde a \in \left(L^\infty(\RR^d)\cap L^r(\RR^d)\right)^{d\times d},\quad  \tilde b \in
  \left(L^\infty(\RR^d)\cap L^s(\RR^d)\right)^d, \\ \hfill \hbox{\rm for some }\quad 1\leq r,s<+\infty, \\
   a^{per},\, \tilde a\in \left(C^{0,\alpha}_{\rm unif}(\RR^d)\right)^{d\times d}, \quad b^{per},\, \tilde b\in \left(C^{0,\alpha}_{\rm unif}(\RR^d)\right)^d \quad\hbox{\rm for some}\quad \alpha>0,
\end{array}
\right.
\end{equation}
We also note that, without loss of generality and because of the specific form of the operator~$-\, a_{ij}\,\partial_{ij} $, we may always assume that $a$ is symmetric. 
We aim to show that the solution to~\eqref{eq:equation-adv-diff-epsilon} may be efficiently approximated using the same ingredients as classical periodic homogenization theory and with the same quality of approximation.
In a series of works~\cite{milan,cras-defauts,cpde-defauts} (see also~\cite{josien,josien-these}), we have studied the same issue for the equation in divergence  form $ -\hbox{\rm div}\, \left(a(x/\varepsilon)\,\nabla u^\varepsilon\right)=f$.  The heart of the matter is the existence  of a corrector function $w_p$, strictly sublinear at infinity (that is,  
$ \displaystyle \frac{ w_{p}(x)}{1+|x|}\quad\buildrel |x|\rightarrow \infty\over\longrightarrow \quad 0$), solution, for each $p\in \RR^d$, to $  -\,\hbox{\rm div}\left(a\,(p+\nabla w_{p})\right)=0$ in $\RR^d$. More precisely, ~$w_p=w_{p,per}+\tilde w_p$ with $w_{p,per}$ the periodic corrector and $\tilde w_p$ solution with $\nabla \tilde w_p\in L^r$ to $ -\,\hbox{\rm div}\left(a\,\nabla \tilde w_{p}\right)= \,\hbox{\rm div}\left(\tilde a\,(p+\nabla w_{p,per})\right)$ in~$\RR^d$.
Such a situation comes in sharp contrast to the general case of homogenization theory where only a sequence of "approximate" correctors is needed to conclude, but where the rate of convergence of the approximation is then unknown. The existence of $\tilde w_p$ above is actually a consequence of the {\it a priori} estimate 
\begin{equation}
  \label{eq:estimee-informelle-div}
  \left\|\nabla u\right\|_{L^q}\,\leq C_q\, \left\|f\right\|_{L^q}, 
\end{equation}
for the exponent $q=r$, and $u$ solution to 
\begin{equation}
  \label{eq:equation-informelle-div}
  -\,\hbox{\rm div}\left(a\,\nabla u\right)= \,\hbox{\rm div}\,f\quad \rm {in }\, \RR^d.
\end{equation}
The case of the equation~\eqref{eq:equation-adv-diff-epsilon} with a vanishing advection field~$b\equiv 0$ has been studied in~\cite{BLL-2017-1}. The corrector associated to this equation identically vanishes, but the issue remains to assess the rate of convergence of the homogenized approximation. This can be achieved proving the estimate
\begin{equation}
  \label{eq:estimee-informelle-non-div-previous}
\left\|D^2\,u\right\|_{\left(L^q(\RR^d)\right)^{d\times d}}\,\leq C_q\, \left\|f\right\|_{L^q(\RR^d)}. 
\end{equation}
for solutions to~ $-\, a_{ij}\,\partial_{ij} u=f$. From that estimate follows the existence of an invariant measure solution to~$-\partial_{ij} (a_{ij}\,m)=0$. Using it to transform equation~\eqref{eq:equation-adv-diff-epsilon}  into an equation in divergence form, we may apply the previous results and conclude. 

\medskip

In all the article, we assume that the dimension satisfies $d\geq 3$. % All our arguments but one, (which is of
% paramount importance, see Remark~\ref{rq:dimension} below) are however valid in the case $d=2$. Of course, the case of dimension one may be dealt with using
% specific tools that we do not wish to develop here. 

\medskip

Our purpose here is to study the general case in~\eqref{eq:equation-adv-diff-epsilon}. Of course, besides~$b\equiv 0$, the  other particular case is when $b_j=\partial_i{a_{ij}}$ in which case equation~\eqref{eq:equation-adv-diff-epsilon} is actually in divergence form~$-\, \partial_i\left(a_{ij}(x/\varepsilon)\,\partial_{j} u^\varepsilon\right) =f$. Otherwise than that, the equation requires a specific treatment.  As in~\cite{BLL-2017-1}, our strategy of proof is based upon establishing an a priori estimate of the type~\eqref{eq:estimee-informelle-non-div-previous}. Because of the presence of the advection field~$b$, a loss in the Lebesgue exponent $q$ will be observed (see  our precise statement in~Proposition~\ref{prop:non-divergence-form} below). Intuitively, and again as in our previous work, the estimate holds true 
because the perturbations~$\tilde a$, and now respectively $\tilde b$,  within the coefficients~$a$ and $b$ respectively, both formally vanish at infinity, while the estimate holds true when $a=a^{per}$, $b=b^{per}$ (using the results of Avellaneda and Lin~\cite{AL1987,AL1989,AL1991}). To the best of our knowledge, it has never been remarked with such a degree of generality that, using an adequate invariant measure, homogenization for the equation~\eqref{eq:equation-adv-diff-epsilon} can be studied and rates can be made precise, simply by transforming the equation into an equation in divergence form. 

\medskip

Our article is organized as follows. We prove in Section~\ref{sec:estimate},
Proposition~\ref{prop:non-divergence-form},  our central estimate. We also explain the loss of integrability we
necessarily observe in comparison to the case of an equation in divergence form or to the case when the advection
field $b$ vanishes. The estimate is then used in Section~\ref{sec:adjoint} to study the adjoint equation
to~\eqref{eq:equation-adv-diff-epsilon}, and prove it admits an invariant measure solution. Various remarks on
possible, very specific cases of coefficients $a$ and $b$ are considered. The invariant measure is in turn employed
in Section~\ref{sec:homog} to transform equation~\eqref{eq:equation-adv-diff-epsilon} in an equation in divergence
form. This allows to apply the results of \cite{cpde-defauts,BLL-2017-1} about the properties of the corrector and
the results of \cite{josien,josien-these} on the approximation of the solution of
\eqref{eq:equation-adv-diff-epsilon} by homogenization theory.

%%%%%%%%%%%%%%%%%%%
\section{The central estimate}
\label{sec:estimate}

We begin by stating and proving our central result (Proposition~\ref{prop:non-divergence-form}) for solutions to
the advection-diffusion equation on the whole space. We will next use the result to prove the existence of an
adequate corrector and conclude the section by some remarks on the optimality of our results. 

\medskip

\begin{proposition}
\label{prop:non-divergence-form}
Assume \eqref{eq:aper+tildea}-\eqref{eq:hyp1} for some $1\leq r< d$ and $\displaystyle 1\leq s<d$. Fix $1\leq q<d$ and set~$\displaystyle  {1\over {q^*}}={1\over q}- {1\over d}$. Then, for all $f\in \left(L^{q^*}\cap L^q\right)(\RR^d)$, there exists~$u\in L^1_{loc}(\RR^d)$ such that $D^2\,u\in L^q(\RR^d)$, solution to
\begin{equation}
\label{eq:equation-adv-diff}
-a_{ij}\partial_{ij}u+b_j\,\partial_j u=f\quad \rm{in}\quad  \RR^d.
\end{equation}
Such a solution is  unique up to the addition of an (at most) affine function. In addition, % the function $u$
                                % belongs to the functional space  {\bf [CLB: Compl\'eter] } and 
there exists a constant $C_q$, independent on $f$ and $u$, and only depending  on $q$, $d$ and the coefficients $a$
and $b$, such that   $u$ satisfies  \begin{equation}
  \label{eq:estimee-informelle-non-div}
\left\|D^2\,u\right\|_{\left(L^{q^*}(\RR^d)\right)^{d\times d}}+ \left\|\nabla u\right\|_{\left(L^{q^*}(\RR^d)\right)^d}\,\leq C_q\, \left\|f\right\|_{\left(L^{q^*}\cap L^q\right)(\RR^d)}.%, \quad {1\over {q^*}}={1\over q}- {1\over d}.
\end{equation}
\end{proposition}

% \begin{remark}
%   {\bf [XB: J'avais note ici qu'il fallait une "remarque sur les $q$ grands". Mais je ne sais plus de quoi il s'agit]}
% \end{remark}

\begin{remark}
In the general case, the above inequality is sharp. However, in the particular case $b=0$, it is not optimal, and this observation is not related to the presence of defects.
  It is already true in the purely periodic case. Indeed, \cite[Proposition 3.1]{BLL-2017-1} gives, in
  the case $b=0$, the estimate, for any $f\in L^q(\RR^d)$,
  \begin{displaymath}
    \left\|D^2\,u\right\|_{\left(L^{q}(\RR^d)\right)^{d\times d}}\,\leq C_q\, \left\|f\right\|_{L^q(\RR^d)},
  \end{displaymath}
where $u$ is the solution of \eqref{eq:equation-adv-diff}. On the other hand, \cite[Theorem B]{AL1991} exactly
states that the result is true (in the periodic case) if and only if the field $b$ vanishes. This will be made
precise in Remarks \ref{rk:perte} and \ref{rk:perte-2} below. Put differently, this loss of decay at infinity is
necessary as soon as a non-trivial transport field $b$ is considered. It is also why we are indeed able to address the case $q=1$ (not covered by Proposition 3.1 of \cite{BLL-2017-1}).  
% Here, we want to be able to deal with the presence of
% $b$. This is why we "loose" some decay at infinity (see Remark~\ref{rk:perte} below on this specific point)
% and are indeed able to address the case $q=1$ (not covered by Proposition 3.1 of \cite{BLL-2017-1}).  
\end{remark}

\noindent{\bf Proof of Proposition~\ref{prop:non-divergence-form}} 

As in~\cite{BLL-2017-1} for the proof of the analogous estimates for the equations in divergence form or the equation~\eqref{eq:equation-adv-diff} with $b\equiv 0$, we argue by continuation. We henceforth fix some~$1\leq q<d$.  We define~$a_t=a^{per}+t\,\tilde a$,~$b_t=b^{per}+t\,\tilde b$  and intend to prove the statements of Proposition~\ref{prop:non-divergence-form} for $t=1$. For this purpose, we introduce the property ${\mathcal P}$ defined by: {\it we say that the coefficients $a$ and $b$, satisfying the assumptions~\eqref{eq:aper+tildea}-\eqref{eq:hyp1} (for some $1\leq r<+\infty$) satisfy ${\mathcal P}$ if the statements of Proposition~\ref{prop:non-divergence-form} hold true for equation~\eqref{eq:equation-informelle-div} with coefficient~$a$ and $b$}. We next define the interval
\begin{equation}
\label{eq:interval}
{\mathcal I}= \left\{t\in [0,1]\,/\, \forall s\in [0,t], \hbox{\rm Property}\,{\mathcal P}\,\hbox{\rm is true
    for}\, a_s\, {\rm and }\,  b_s\right\}.
\end{equation}
We intend to successively prove that ${\mathcal I}$ is not empty, open and closed (both notions being understood relatively to the closed interval $[0,1]$), which will show that ${\mathcal I}=[0,1]$, and thus the result claimed.

\medskip

\noindent {\it Step 1: $0\in {\mathcal I}$.} To start with, we show that $0\in {\mathcal I}$. In the particular case when $b^{per}\equiv 0$ (a case considered in~\cite[Proposition 3.1]{BLL-2017-1}), the fact that~$0\in {\mathcal I}$ is shown to be a consequence of the results of~\cite[Theorem B]{AL1991}. Indeed, the adjoint equation~\eqref{eq:general-equation-adjoint} associated to~\eqref{eq:equation-adv-diff}, which  reads as~$-\partial_{ij}(a^{per}_{ij}\,m_{per})=0$, admits (see e.g.~\cite{blp}), a unique nonnegative periodic solution $m_{per}$ that is normalized, regular and  bounded away from zero. Multiplying~\eqref{eq:equation-adv-diff} by $m_{per}$, we may  write this equation in the divergence form 
\begin{equation}
\label{eq:forme-div-periodique}
-\div\left(\mathcal A^{per} \nabla u\right)=m_{per} f,
\end{equation}
with 
\begin{equation}
\label{eq:mathcal-A-per}
\mathcal A^{per} = m_{per}\,a^{per}- \mathcal B^{per},
\end{equation}
 $\mathcal B^{per}$ the skew-symmetric matrix defined by $\div \left(\mathcal B^{per}\right)=\div(m_{per}\,a^{per})$, and where 
\begin{equation}
\label{eq:div-A-per-nulle}
\div\mathcal A^{per}=\div(m_{per}\,a^{per})-\div\mathcal B^{per}=0.
\end{equation}
A proof of the existence (and uniqueness) of $\mathcal B^{per}$ may be found in \cite[Chapter 1]{jikov}.
We thus deduce from~\cite[Theorem B]{AL1991} that  
\begin{equation}
  \label{eq:estim_q}
 \left\|D^2\,u\right\|_{\left(L^q(\RR^d)\right)^{d\times
    d}}\,\leq C\, \left\|f\right\|_{L^q(\RR^d)}.
\end{equation}
This inequality is actually valid for any $q>1$, hence it holds
also  for $q^*$: 
\begin{displaymath}
 \left\|D^2\,u\right\|_{\left(L^{q^*}(\RR^d)\right)^{d\times
    d}}\,\leq C\, \left\|f\right\|_{L^{q^*}(\RR^d)}.
\end{displaymath}
Using \eqref{eq:estim_q} and Gagliardo-Nirenberg-Sobolev inequality (see for instance \cite[Section 5.6.1,
Theorem 1]{evans}), we infer $\|\nabla u\|_{\left(L^{q^*}(\RR^d)\right)^d} \leq C\|f\|_{L^q(\RR^d)}.$
The local integrability $u\in L^1_{loc}(\RR^d)$ is obtained by elliptic regularity  using $f\in L^1_{loc}(\RR^d)$ and the H\"older regularity of the coefficient~$a^{per}$ stated in~\eqref{eq:hyp1}. This property immediately carries over to  all the other cases we henceforth consider as soon we know there is a solution.

We next insert a non vanishing advection field $b^{per}$. Unless $b^{per}_j=\partial_i\,a^{per}_{ij}$ (and the
equation is then in divergence form), we have to work more. We still have,  as above again because of the classical
results exposed in~\cite{blp}, the existence of an invariant measure, this time solution
to~~$-\partial_{j}(\partial_i(m_{per}\,a^{per}_{ij})+\,m_{per}b^{per}_j)\,=0$, with all the suitable
properties. This allows again to write the original equation in the divergence
form~\eqref{eq:forme-div-periodique}, but this time, \eqref{eq:div-A-per-nulle} is not satisfied and we cannot
apply~\cite[Theorem B]{AL1991}.  However, since~\eqref{eq:forme-div-periodique} holds, and since the matrix-valued
coefficient is periodic and regular (because of~\eqref{eq:hyp1}), we know that the Green function $G^{per}(x,y)$
associated to the operator~$-\div\left(\mathcal A^{per} \nabla \, .\right)$ satisfies, for all $x,y\in\RR^d$, (see
\cite[Theorem 1.1]{GW} and \cite[Proposition 2]{anantha})
% \begin{equation}
% \label{eq:marcinkiewitz}
% \nabla G^{per} (x,.)\,\in\, L^{d/(d-1),\infty}(\RR^d)
% \end{equation}
\begin{equation}
\label{eq:marcinkiewitz}
\left|\nabla G^{per} (x,y)\right| \leq \frac C {|x-y|^{d-1}}\ .%\,\in\, L^{d/(d-1),\infty}(\RR^d)
\end{equation}
Hence,
\begin{displaymath}
  \nabla u(x) = \int_{\RR^d} \nabla_x G^{per}(x,y)\ m_{per}(y)f(y)dy
\end{displaymath}
satisfies
\begin{displaymath}
  |\nabla u(x)|\leq  \left\|m_{per}\right\|_{L^{\infty}(\RR^d)}\, \int_{\RR^d} \frac C {|x-y|^{d-1}} |f(y)|dy.
\end{displaymath}
Now, the O'Neil-Young inequality \cite{oneil,yap} states that $\forall f\in L^{p_1,q_1}(\RR^d),$ $\forall g\in L^{p_2,q_2}(\RR^d),$
\begin{equation}
  \label{eq:young-oneil}
  \|f*g\|_{L^{\sigma,\theta}(\RR^d)} \leq C
  \|f\|_{L^{p_1,q_1}(\RR^d)} \,  \|g\|_{L^{p_2,q_2}(\RR^d)},
\end{equation}
where $\frac 1 {p_1} + \frac 1 {p_2} = 1+ \frac 1 \sigma$ and $\frac 1 {q_1} + \frac 1 {q_2}\geq \frac 1 \theta$, $1\leq p_i
\leq \infty$, $1\leq q_i \leq \infty$ (except for the case $(p_i=1,q_i=\infty)$)and $L^{p,q}$
denotes the Lorentz space of exponent $(p,q)$ (see \cite{grafakos,lorentz}.) The constant $C$
in \eqref{eq:young-oneil} does not depend on $f$ and $g$. It is easily proved that $|x|^{-(d-1)} \in
L^{d/(d-1),\infty}(\RR^d)$, hence, since $f\in L^q(\RR^d)$,
\begin{equation}
\label{eq:convol0}
\left\|\nabla u\right\|_{\left(L^{q^*,\theta}(\RR^d)\right)^d}\leq C\,
\left\|m_{per}\right\|_{L^{\infty}(\RR^d)}\,\left\|f\right\|_{L^{q}(\RR^d)} \, \left\| \frac 1 {|x|^{d-1}}\right\|_{L^{d/(d-1),\infty}(\RR^d)},
\end{equation}
provided $\frac 1 \theta \leq \frac 1 q$. Since $\frac 1 {q^*} = \frac 1 q - \frac 1 d$, $\theta=q^*$ is allowed in
\eqref{eq:young-oneil}. Therefore,
\begin{equation}\label{eq:convol}
\left\|\nabla u\right\|_{\left(L^{q^*}(\RR^d)\right)^d}\leq C\,
\left\|m_{per}\right\|_{L^{\infty}(\RR^d)}\,\left\|f\right\|_{L^{q}(\RR^d)}.
\end{equation}
% Since the right-hand side of \eqref{eq:marcinkiewitz} is in the Marcinkiewitz space $L^{d/(d-1),\infty}(\RR^d)$, 

% Consequently, applying the H\"older inequality and using $m_{per}\,f\in 
% L^q(\RR^d)$, 
% % \begin{equation}
% % \label{eq:convol}
% % \left\|\nabla u=\int \nabla G^{per} (.,y)\,(m_{per}(y)\,f(y))\,dy\right\|_{\left(L^{q^*}(\RR^d)\right)^d}\leq C\, \left\|m_{per}\right\|_{L^{\infty}(\RR^d)}\,\left\|f\right\|_{L^{q}(\RR^d)}
% % \end{equation}
% since $\displaystyle  1+ {1\over {q^*}}={{(d-1)}\over d}+ {1\over q}$.
We next rewrite $-a^{per}_{ij}\partial_{ij}u+b^{per}_j\,\partial_j u=f$ as~$-a^{per}_{ij}\partial_{ij}u=f-b^{per}_j\,\partial_j u$.  In the right-hand side of the latter equation, we note that
\begin{equation}
\label{eq:rhs}
\left\|f-b^{per}_j\,\partial_j u\right\|_{L^{q^*}(\RR^d)}\leq \, \left\|f\right\|_{L^{q^*}(\RR^d)}\,+\,\left\|b^{per}\right\|_{L^{\infty}(\RR^d)}\,\left\|\nabla u\right\|_{\left(L^{q^*}(\RR^d)\right)^d}
\end{equation}
%
% Ici on a besoin de $f\in L^{q^*}$ en plus de $f\in L^q$, 
% et qui est de toute façon une condition n\'ecessaire! 
% Pour la suite (existence de $m$), \`a verifier
%
We may therefore apply~\cite[Theorem B]{AL1991}: inserting~\eqref{eq:convol} into~\eqref{eq:rhs}, we obtain~\eqref{eq:estimee-informelle-non-div} in the specific case of periodic coefficients.

\medskip

\noindent {\it Step 2: ${\mathcal I}$ is open.} The fact  that ${\mathcal I}$ is open (relatively to the interval~$[0,1]$) is a straightforward consequence of the Banach fixed point Theorem. We solve, for $f\in L^q(\RR^d)$ fixed and $\varepsilon>0$ presumably small,
$$-((a_t)_{ij}+ \varepsilon\,\tilde a_{ij})\partial_{ij}u+((b_t)_j+ \varepsilon\,\tilde b_{j})\,\partial_j u=f\quad \rm {in}\quad \RR^d,$$
using the iterations~$u^0=0$ and, for all~$n\in\NN$,  
$$-(a_t)_{ij}\partial_{ij}u^{n+1}+(b_t)_j\,\partial_j u^{n+1}=f \,+\,\varepsilon\,\left(\tilde a_{ij}\partial_{ij}u^{n}-\,\tilde b_{j}\,\partial_j u^{n}\right)
$$
The point is to prove that the right-hand side belongs to  $\left(L^{q^*}\cap L^q\right)(\RR^d)$. By assumption,
$f\in \left(L^{q^*}\cap L^q\right)(\RR^d)$. We also have, by the H\"older inequality,  (i) because $r\leq d$, $\tilde a \in \left(\left(
    L^d\cap L^{\infty}\right)(\RR^d)\right)^{d\times d}$ and, by inductive hypothesis,
$D^2u^{n}\in\left(L^{q^*}(\RR^d)\right)^{d\times d}$, thus $\tilde a:D^2u^{n}\in \left(L^q\cap
  L^{q^*}\right)(\RR^d)$, (ii) because $s\leq d$, $\tilde b \in \left(L^{d}\cap L^\infty(\RR^d)\right)^{d}$ and, by inductive
hypothesis and the Sobolev embedding Theorem, $\nabla u^{n}\in\left(L^{q^{*}}(\RR^d)\right)^{d}$ (we recall that $\displaystyle {1\over {q^{*}}}= {1\over {q}} - {1\over d}$,) thus $\tilde b\,.\,\nabla u^{n}\in \left(L^q\cap L^{q^*}\right)(\RR^d)$. 
 By induction, the iterate~$u^{n+1}$ is thus well defined (up to an irrelevant, at most affine, function)
 with~$D^2u^{n+1}\in\left(L^{q^*}(\RR^d)\right)^{d\times d}$ and $\nabla u^{n+1} \in \left(L^{q^*}(\RR^d)\right)^{d}$, precisely applying Property ${\mathcal P}$ for the coefficients~$a_t$, $b_t$. Also because of that property, we have, for $\varepsilon$ sufficiently small, a geometric convergence of the series $\sum_n(u^{n+1}-u^{n})$. Existence of the solution $u$ follows. The uniqueness of a solution (again up to the addition of an irrelevant at most affine function) is proven similarly.
 
\medskip

\noindent {\it Step 3: ${\mathcal I}$ is closed.} We now show, and this is the key point of the proof,  that
${\mathcal I}$ is closed. We assume that $t_n\in {\mathcal I}$, $t_n\leq t$, $t_n\longrightarrow t$ as
$n\longrightarrow +\infty$. For all $n\in\NN$, we know that, for any $f\in \left(L^{q^*}\cap L^q\right)(\RR^d)$, we
have a solution (unique to the addition of an irrelevant function) $u^n$ with~$D^2 u^n\in
\left(L^{q^*}(\RR^d)\right)^{d\times d}$ and $\nabla u \in \left(L^{q^*}(\RR^d)\right)^d$ of the equation 
$$-(a_{t_n})_{ij}\partial_{ij}u^n+(b_{t_n})_j\,\partial_j u^n=f\quad \rm {in }\quad \RR^d,
$$
and that this solution satisfies 
\begin{displaymath}
 \left\|D^2 u^n\right\|_{\left(L^{q^*}(\RR^d)\right)^{d\times d}} +\left\|\nabla u^n\right\|_{\left(L^{q^*}(\RR^d)\right)^d} \leq C_n\,\left\|f\right\|_{\left(L^{q^*}\cap L^q\right)(\RR^d)},
\end{displaymath}
for a constant $C_n$ depending on $n$ but not on~$f$ nor on~$u^n$. We want to show the same properties for $t$. 

\medskip

We first conclude temporarily admitting that  the constants $C_n$ are bounded uniformly in~$n$. Next, we will prove this is indeed the case. For $f\in \left(L^{q}(\RR^d)\right)^{d}$ fixed, we consider the sequence of solutions~$u^n$ to 
$$-(a_{t_n})_{ij}\partial_{ij}u^n+(b_{t_n})_j\,\partial_j u^n=f\quad {\rm in }\quad \RR^d,
$$
which we may write as
$$-(a_t)_{ij}\partial_{ij}u^n+(b_t)_j\,\partial_j u^n=f+(t-t_n)\left(-\tilde a_{ij}\partial_{ij}u^n+\tilde b_{j}\,\partial_j u^n\right)\quad {\rm in }\quad \RR^d,
$$
Since $t_n\in{\mathcal I}$ for all $n\in\NN$ and the constants $C_n$ are uniformly bounded, we know that the
sequences $D^2 u^n$ and $\nabla u^n$ are bounded in~$\left(L^{q^*}(\RR^d)\right)^{d\times d}$ and
in~$\left(L^{q^{*}}(\RR^d)\right)^{d}$, respectively. We may pass  to the weak limit in the above equation and find
a solution~$u$ to $-(a_t)_{ij}\partial_{ij}u+(b_t)_j\,\partial_j u=f$. The solution also satisfies the estimate
(because the sequence $C_n$ is bounded and because the norm is weakly lower semi continuous). 

\medskip

In order to prove that the constants $C_n$ are indeed bounded uniformly in~$n$, we argue by contradiction. We assume we have ~$f^n\in \left(L^{q}(\RR^d)\right)^{d}$ and $u^n$ with~$D^2 u^n\in \left(L^{q}(\RR^d)\right)^{d\times d}$, such that 
\begin{equation}
\label{eq:contradic1-nondiv}
-(a_{t_n})_{ij}\partial_{ij}u^n+(b_{t_n})_j\,\partial_j u^n= \,f^n\quad {\rm in }\quad \RR^d,
\end{equation}
\begin{equation}
\label{eq:contradic2-nondiv}
\left\|f^n\right\|_{\left(L^{q^*}\cap L^q\right)(\RR^d)}\buildrel n\longrightarrow +\infty\over \longrightarrow 0,
\end{equation}
\begin{equation}
\label{eq:contradic3-nondiv}
\left\|\nabla u^{n}\right\|_{\left(L^{q^*}(\RR^d)\right)^d}+\left\|D^2 u^{n}\right\|_{\left(L^{q^*}(\RR^d)\right)^{d\times d}}=1,\quad\hbox{\rm for all}\,n\in \NN.
\end{equation}
To start with, we rewrite~\eqref{eq:contradic1-nondiv} as
$$
-(a_t)_{ij}\partial_{ij}u^n+(b_t)_j\,\partial_j u^n= \,f^n+ (t_n-t)\,\tilde a_{ij}\,\partial_{ij}u^n+(t_n-t)\,\tilde b_{j}\,\partial_j u^n,
$$
where,  as $n\longrightarrow 0$, the rightmost two terms vanish in $\left(L^{q^*}\cap L^q\right)(\RR^d)$ using the bound~\eqref{eq:contradic3-nondiv} and the  same argument as above for the openness of ${\mathcal I}$. Therefore, without loss of generality, we may change the definition of $f^n$ and  replace~\eqref{eq:contradic1-nondiv} by
\begin{equation}
\label{eq:contradic1-bis-nondiv}
-(a_t)_{ij}\partial_{ij}u^n+(b_t)_j\,\partial_j u^n= \,f^n\quad {\rm in }\, \RR^d,
\end{equation}
In the spirit of the method of concentration-compactness \cite{PLL-cc}, we now claim that the sequence $u^{n}$ satisfies 
\begin{multline}
\label{eq:cc-1-nondiv}
\exists\, \eta>0, \quad \exists\, 0<R<+\infty,\quad \forall \,n\in\NN,\quad \\ \left\|D^2
  u^{n}\right\|_{\left(L^{q^*}(B_R)\right)^{d\times d}}+ \left\|\nabla u^{n}\right\|_{\left(L^{q^*}(B_R)\right)^d}\geq \eta>0,
\end{multline}
where $B_R$ of course denotes the ball of radius~$R$ centered at the origin.  We again argue by contradiction  and assume that, contrary to~\eqref{eq:cc-1-nondiv},
\begin{equation}
\label{eq:cc-1-not}
\forall\, 0<R<+\infty,\quad  \left\|D^2 u^{n}\right\|_{\left(L^{q^*}(B_R)\right)^{d\times d}}+ \left\|\nabla u^{n}\right\|_{\left(L^{q^*}(B_R)\right)^d}\buildrel n\longrightarrow +\infty\over \longrightarrow 0.
\end{equation}
Since both $\tilde a$ and $\tilde b$ satisfy the properties in~\eqref{eq:hyp1}, they vanish at infinity and thus, for any~$\delta>0$,  we may find some sufficiently large radius~$R$ such that 
\begin{equation}
\label{eq:cc-2}
\left\|\tilde a\right\|_{\left(\left(L^d\cap L^\infty\right)(B_R^c)\right)^{d\times d}}\leq \delta , \quad \left\|\tilde b\right\|_{\left(\left(L^{d}\cap L^\infty\right)(B_R^c)\right)^{d}}\leq \delta ,
\end{equation}
where~$B_R^c$ denotes the complement set of the ball $B_R$. We then estimate \begin{eqnarray}\label{eq:5}
\left\|\tilde a\,D^2 u^n\right\|^{q}_{L^{q}(\RR^d)}&=&\left\|\tilde a\,D^2 u^n\right\|^{q}_{L^{q}(B_R)}+\left\|\tilde a\,D^2 u^n\right\|^{q}_{L^{q}(B_R^c)}\nonumber\\
&\leq&\left\|\tilde a\right\|^q_{\left(L^d(\RR^d)\right)^{d\times d}}\,\left\|D^2 u^n\right\|^{q^*}_{\left(L^{q^*}(B_R)\right)^{d\times d}}\nonumber\\
&&
+\left\|\tilde a\right\|^q_{\left(L^d(B_R^c)\right)^{d\times d}}\,\left\|D^2 u^n\right\|^{q^*}_{\left(L^{q^*}(\RR^d)\right)^{d\times d}}\nonumber\\
&\leq&\left\|\tilde a\right\|^{q}_{\left(L^d(\RR^d)\right)^{d\times d}}\,\left\|D^2u^n\right\|^q_{\left(L^{q^*}(B_R)\right)^{d\times d}}%\nonumber\\
%&&
+\delta ,
\end{eqnarray}
using~\eqref{eq:contradic3-nondiv} and~\eqref{eq:cc-2} for the latter majoration. Given that~\eqref{eq:cc-1-not}
implies that the first term in the right hand side of \eqref{eq:5} vanishes, and since $\delta$ is arbitrary, this
shows that~$\tilde a\,D^2 u^n\to 0$ in $L^{q}(\RR^d)$. By the exact same argument, this time using the $L^\infty$ estimate of $\tilde a$ in~\eqref{eq:cc-2}, we likewise obtain that $\tilde a\,D^2 u^n$ vanishes in $L^{q^*}(\RR^d)$. Therefore
\begin{equation}
\label{eq:a-vanish}
 \left\|\tilde a\,D^2 u^{n}\right\|_{\left(L^{q^*}\cap L^q\right)(\RR^d)}\buildrel n\longrightarrow +\infty\over \longrightarrow 0.
\end{equation}
We address the first order term similarly, getting
\begin{equation}
\label{eq:b-vanish}
 \left\|\tilde b\,\nabla u^{n}\right\|_{\left(L^{q^*}\cap L^q\right)(\RR^d)}\buildrel n\longrightarrow +\infty\over \longrightarrow 0.
\end{equation}
% On the one hand because of the estimates on $\tilde b$
                                % in~\eqref{eq:cc-2}, and on the other hand because~\eqref{eq:cc-1-not} and the Sobolev embedding Theorem implies that 
% \begin{equation}
% \label{eq:cc-1-not-grad}
% \forall\, 0<R<+\infty,\quad  \left\|\nabla u^{n}\right\|_{\left(L^{q^{**}}(B_R)\right)^{d}}\buildrel n\longrightarrow +\infty\over \longrightarrow 0
% \end{equation}
%  {\bf [CLB : V\'erifier que~\eqref{eq:cc-1-not} implique bien ~\eqref{eq:cc-1-not-grad}, modulo un bon choix de "normalisation" de~$u^{n}$  !!]} where, we recall, $\displaystyle {1\over {q^{**}}}= {1\over {q^*}} - {1\over d}$, we have
% \begin{equation}
% \label{eq:b-vanish}
%  \left\|\tilde b\,\nabla u^{n}\right\|_{\left(L^{q^*}\cap L^q\right)(\RR^d)}\buildrel n\longrightarrow +\infty\over \longrightarrow 0.
% \end{equation}
We next notice that~\eqref{eq:contradic1-bis-nondiv} also reads as 
 $$
-a^{per}_{ij}\partial_{ij}u^n+b^{per}_{j}\,\partial_j u^n= \,f^n+ t\,\tilde a_{ij}\,\partial_{ij}u^n-t\,\tilde b_{j}\,\partial_j u^n,
$$
and use~\eqref{eq:contradic2-nondiv},~\eqref{eq:a-vanish} and~\eqref{eq:b-vanish} to estimate its right-hand
side. In view of the estimate~\eqref{eq:estimee-informelle-non-div} which, as mentioned above, holds for the
operator for periodic coefficients, this implies that $D^2 u^n$ and $\nabla u^n$ (strongly) converge to zero
in~$\left(L^{q^*}(\RR^d)\right)^{d\times d}$ and $\left(L^{q^*}(\RR^d)\right)^d$, respectively. This evidently contradicts~\eqref{eq:contradic3-nondiv} . We
therefore have established~\eqref{eq:cc-1-nondiv}.

\medskip

We are now in position to finally reach a contradiction.
Because of the bound~\eqref{eq:contradic3-nondiv}, we may claim that, up to an extraction, $D^2 u^n$ weakly
converges in $\left(L^{q^*}(\RR^d)\right)^{d\times d}$, to some $D^2 u$. This convergence is actually  strong in
$\left(L^{q^*}_{loc}(\RR^d)\right)^{d\times d}$. This is proven combining  Sobolev compact embeddings and
estimates for general elliptic operators (see e.g.~\cite[Theorem 7.3]{Giaquinta-Martinazzi}).  Passing to the weak limit
in~\eqref{eq:contradic1-bis-nondiv}, we obtain $-\,(a_t)_{ij}\partial_{ij} u+(b_t)_j\partial_{j} u= 0$ for $u$ that
does not  identically vanish. This is a contradiction with the uniqueness we prove below. 

\medskip

There remains to prove uniqueness. We thus consider a solution $u$ to \eqref{eq:equation-adv-diff} with $f=0$,
$D^2 u \in L^{q^*}(\RR^d)^{d\times d}$ and $\nabla u \in L^{q^*}(\RR^d)^d$. 

We first consider the case $q<d/2$,
i.e $q^*<d$. In such a case, the Gagliardo-Nirenberg-Sobolev inequality implies that, up to the addition of a
constant, $u\in L^{q^{**}}(\RR^d)$. Moreover, using elliptic regularity \cite[Theorem 9.11]{GT}, one easily proves
that
\begin{equation}\label{eq:8}
  \sup_{x\in\RR^d} \|u\|_{W^{2,q^{**}}(B_1(x))} <+\infty,
\end{equation}
where we recall that $\frac 1 {q^{**}} = \frac 1 q - \frac 2 d$. If $q^{**}>d$, we apply Morrey's theorem, proving that
$u$ is H\"older continuous. If not, we repeat the above argument, obtaining \eqref{eq:8}
with $q^{***}$, and so on, that is, for any integer $n$:
\begin{displaymath}
  \sup_{x\in\RR^d} \|u\|_{W^{2,q_n}(B_1(x))} <+\infty, \quad \frac 1 {q_n} = \frac 1 q - \frac n d, \quad \text{as
    long as}\quad n < \frac d q.
\end{displaymath}
For $m$ such that $m< \frac d q < m +1$ (if $\frac d q \in \NN,$ slightly decrease $q$ such that it is no longer
the case; this is possible because the estimate is local), we have $q_m >d$, hence, by Morrey's
theorem, $u\in C^{0,\alpha}_{\rm unif}(\RR^d)$, for some $\alpha>0$. This and $u\in L^{q^{**}}(\RR^d)$ imply that, for any
$\delta >0$, there exists $R>0$ such that
\begin{displaymath}
  \sup_{|x|>R} |u(x)|\leq \delta.
\end{displaymath}
Applying the maximum principle on the ball $B_R$, we thus have $|u|\leq \delta$ in $\RR^d$. Since this is valid for
any $\delta>0$, we conclude that $u=0$. 

In order to address the case $d/2 \leq q<d $, we write \eqref{eq:equation-adv-diff} as
\begin{equation}\label{eq:6}
  -a_{ij}^{per} \partial_{ij}u + b_j^{per} \partial_j u =  t\left( \tilde a_{ij}\partial_{ij} u - \tilde b_j\partial_j u\right) 
\end{equation}
Here again, the fact that $\tilde a \in \left(
    L^r\cap L^\infty(\RR^d)\right)^{d\times d}$ and $D^2 u \in L^{q^*}(\RR^d)^{d\times d}$ implies that
$\tilde a_{ij} \partial_{ij} u \in L^{r_1}\cap L^{q^*}(\RR^d)$, with $\frac 1 {r_1} = \frac 1 r + \frac 1 {q^*} = \frac 1 r +
\frac 1 {q} - \frac 1 d.$ Similarly, $\tilde b_j\partial_j u\in L^{s_1}\cap L^{q^*}(\RR^d)$, with $\frac 1 {s_1} =\frac 1 s +
\frac 1 {q} - \frac 1 d.$ Since $r,s<d$, we have $r_1,s_1<q$. Applying step 1 of the present proof,
we conclude that
\begin{equation}\label{eq:9}
  D^2u \in \left(L^{\max(r_1,s_1)^*}\cap L^{q^*}(\RR^d)\right)^{d\times d}, \quad   \nabla u \in \left(L^{\max(r_1,s_1)^*}\cap L^{q^*}(\RR^d)\right)^{d}.
\end{equation}
Repeating this argument, \eqref{eq:9} is also valid for $s_n$ and $r_n$ defined by 
\begin{displaymath}
  \frac 1 {s_n} = \frac 1 q + \frac n s - \frac n d,\quad \frac 1 {r_n} = \frac 1 q + \frac n r - \frac n d.
\end{displaymath}
Hence, for $n$ sufficiently large, we have $\max(r_n,s_n)<d/2$, and we may apply the argument of the case $q<d/2$.

\medskip

We reach a final contradiction. This shows that  ${\mathcal I}$ is closed. As it is also open and non empty, it is equal to~$[0,1]$ and this concludes the proof of Proposition~\ref{prop:non-divergence-form}. \hfill $\diamondsuit$

% \bigskip
% \begin{remark}\label{rq:dimension}
%   In the above proof, the uniqueness argument starting with \eqref{eq:8} is only valid for $d\geq 3$. Indeed, in
%   the case $d=2$, the assumption $q<d/2$ is in contradiction with $q>1$. We do not know how to adapt it to the case
%   $d=2$. Les pistes explorees:
%   \begin{enumerate}
%   \item Utiliser Cordes comme dans le cas $b=0$. Ici, ca donne
%     \[ \left(b\cdot\nabla u\right)^2 = \left(a_{ij}\partial_{ij} u \right)^2 \geq C_0 |D^2u|^2 + C_1 \det(D^2 u).\]
% Cette inegalite est ponctuelle, mais pour faire disparaitre le dernier terme, dont on ne connait pas le signe, il
% faut integrer sur des grandes boules. Ensuite, meme si on se debarasse de ce terme, on ne sait pas vraiment
% utiliser l'inegalite qui en decoule. 
% \item Utiliser l'inegalite de Harnack (elle est demontree pour un tel operateur dans un papier de Trudinger de
%   1980). Le probleme : la constante de Harnack depend du domaine, contrairement au cas $b=0$. On ne peut donc pas
%   appliquer la strategie du papier de Moser qui montre qu'une solution de l'equation homogene croit au moins comme
%   une puissance a l'infini, si elle n'est pas constante.
%   \end{enumerate}
% \end{remark}
\bigskip

\begin{remark}
  It is clear from the above proof that Proposition~\ref{prop:non-divergence-form} is also valid in the case $r=d$ and/or $s=d$ if we assume in
  addition that $q<\frac d 2$. The only stage where $s,r<d$ is used is to prove uniqueness in the case $q\geq
  \frac d 2$.
\end{remark}

We now use Proposition~\ref{prop:non-divergence-form} to prove the existence of a corrector for our problem. We
first recall the following facts for the periodic case (see e.g. \cite{blp}). There exists a unique positive
measure, bounded away from zero, with normalized periodic average $\left\langle m_{per}\right\rangle =1$, that solves
\begin{equation}
\label{eq:adjoint-per}
-\partial_i\left(a^{per}_{ij}\partial_{j}m_{per}+b^{per}_j\,m_{per}\right)=0\quad \rm{in }\, \RR^d.
\end{equation}
If the condition 
\begin{equation}
\label{eq:mper-bper}
\left\langle m_{per} b^{per}\right\rangle=0
\end{equation}
holds true, then, for all $p\in\RR^d$, there exists a periodic corrector function, with normalized average
$\left\langle w_{p,per}\right\rangle =0$, solution to 
\begin{equation}
\label{eq:corrector-per}
-a^{per}_{ij}\partial_{ij}w_{p,per}+b^{per}_j\,\partial_j w_{p,per}=-b^{per}.p\quad \rm{in }\, \RR^d.
\end{equation}
By elliptic regularity (see for instance \cite[Theorem 9.11]{GT}), this function satisfies $w_{p,per}\in L^\infty(\RR^d)$, $\nabla\,w_{p,per}\in \left(L^\infty(\RR^d)\right)^d$ and $D^2 w_{p,per}\in \left(L^\infty(\RR^d)\right)^{d\times d}$. 

\medskip

\begin{corollary}
\label{cor:corrector}
As in Proposition~\ref{prop:non-divergence-form}, we assume \eqref{eq:aper+tildea}-\eqref{eq:hyp1} for some $1\leq r<d$ and $\displaystyle 1\leq s< d$.  We additionally  assume % that $r<d$ 
% % donc $\max(r,s)<d$. 
% and
the condition~\eqref{eq:mper-bper}. Then, for all $p\in\RR^d$, there exists a corrector function, solution to 
\begin{equation}
\label{eq:corrector}
-a_{ij}\partial_{ij}w_p+b_j\,\partial_j w_p=-b.p\quad \rm{in }\, \RR^d.
\end{equation}
Such a solution is  unique up to the addition of an (at most) affine function. It reads as \begin{equation}
\label{eq:corrector-decomp}
w_p=w_{p,per}+\tilde w_p
\end{equation}
where $w_{p,per}$  is the periodic corrector solution to \eqref{eq:corrector-per} with normalized average $\langle
w_{p,per}\rangle=0$, and where (again up to the addition of an at most affine function) $\tilde w_p\in
L^1_{loc}(\RR^d)$, $\nabla\tilde w_p \in \left(L^{q^*}(\RR^d)\right)^d$, $D^2\,\tilde w_p\in
\left(L^{q^*}(\RR^d)\right)^{d\times d}$ for ${1\over {q^*}}={1\over {\max(r,s)}}- {1\over d}$. In particular, the corrector $w_p$ is thus strictly sub-linear at infinity. % Also, when $d\geq 3$ and $\max(r,s)<d/2$,  
% \begin{equation}
% \label{eq:tilde-w-Lq}
% \nabla\,\tilde w_p\in L^{q^{**}}(\RR^d)\quad\hbox{\rm for}\quad {1\over {q^{**}}}=\min\left({1\over r},{1\over s} \right)- {2\over d},
% \end{equation}
% and the corrector $\tilde w_p$ is strictly sublinear at infinity.
\end{corollary}

\medskip

\noindent{\bf Proof of Corollary~\ref{cor:corrector}}
Using~\eqref{eq:corrector-per}, we notice that~\eqref{eq:corrector} also reads as
\begin{equation}
\label{eq:corrector-bis}
-a_{ij}\partial_{ij}\tilde w_p+b_j\,\partial_j \tilde w_p=-\tilde b.p+\tilde a_{ij}\partial_{ij}w_{p,per}-\tilde b_j\,\partial_j w_{p,per}\quad \rm{in }\, \RR^d.
\end{equation}
Given the properties of boundedness of~$w_{p,per}$ and its first and second derivatives and our assumptions~\eqref{eq:aper+tildea}-\eqref{eq:hyp1}, the right-hand side of~\eqref{eq:corrector-bis} belongs to $\left(L^{\max(r,s)}\cap L^\infty\right)(\RR^d)$. Since we have assumed $\max(r,s)<d$, we may apply Proposition~\ref{prop:non-divergence-form} for the exponent $q=\max(r,s)$ and we obtain the results stated in Corollary~\ref{cor:corrector}.  \hfill $\diamondsuit$

\medskip

We conclude this section with a series of remarks on our assumptions and results of Proposition~\ref{prop:non-divergence-form} and Corollary~\ref{cor:corrector}.

\medskip

\begin{remark}
% {\bf [CLB: COMMENTER LA PERTE D'EXPOSANT (on DOIT perdre en $r$ sinon on pourrait consid\'erer l'op\'erateur
%   rescal\'e, et on fait exploser l'estim\'ee avec un scaling en $\varepsilon$)), l'argument est identique \`a celui
%   fait dans les Notes XB du 17/03/2017: on montre ainsi que \underbar{si} une in\'egalit\'e du
%   type~\eqref{eq:estimee-informelle-non-div} pour un certain $L^\alpha$ est vraie, \underbar{alors}
%   n\'ecessairement $\alpha={q^*}$] 

% [XB: J'ai ecrit ce qui suit, mais je ne suis pas sur que ce soit ce que tu avais
%   en tete]}
One should not be surprised by the fact that, in the left-hand side of \eqref{eq:estimee-informelle-non-div}, the first derivative
$\nabla u$ and the second derivative $D^2u$ share the \emph{same} integration exponent. One could think that,
because of Gagliardo-Nirenberg-Sobolev inequality, the exponent of the first derivative and that of the second
derivative are related by $\frac 1 {q^*} = \frac 1 q - \frac 1 d$. Because of the \emph{structure} of the
differential operator, it is indeed possible to have the same exponent. To illustrate this idea, we consider the
simple example where $a = \mathbf{Id}$, and $b_i(x) = b_0(|x|)\frac{x_i}{|x|}$. We assume
$r\mapsto b_0(r)$ to be smooth and vanish at $0$ in order to have $b\in C^{0,\alpha}_{\rm unif}(\RR^d)$, and $b_0(r)=1$
for $r\geq 1$. We also assume that an estimate of the form $\|D^2u\|_{L^\beta} \leq C\|f\|_{L^q}$ holds for the solution of
\eqref{eq:equation-adv-diff}, for some exponent $\beta$. If the right-hand side $f$ is radially symmetric, so is the solution $u$, and the equation reads
\begin{equation}\label{eq:10}
  -\frac{d^2 u}{dr^2} - \frac{d-1}r \frac{du}{dr} + b_0(r)\frac{du}{dr} = f(r).
\end{equation}
Since $b_0(r) + \frac{d-1}r \to 1$ as $r\to+\infty$, the estimate $\|D^2u\|_{L^\beta} \leq C\|f\|_{L^q}$, together
with equation \eqref{eq:10}, imply $\|\nabla
u\|_{L^\beta}\leq C\|f\|_{L^q}$. 
\end{remark}

\begin{remark}\label{rk:perte} One should not be surprised either by the fact that the exponent of the rignt-hand
  side of \eqref{eq:estimee-informelle-non-div} is equal to $q^*$, which is larger than $q$. Indeed, this is
  already a necessary condition in the periodic case: if we assume that, for equation \eqref{eq:equation-adv-diff}, an estimate of the type
\begin{displaymath}
\left\|D^2\,u\right\|_{\left(L^{\beta}(\RR^d)\right)^{d\times d}}+ \left\|\nabla u\right\|_{\left(L^{\beta}(\RR^d)\right)^d}\,\leq C_q\, \left\|f\right\|_{\left(L^{\beta}\cap L^q\right)(\RR^d)},
\end{displaymath}
holds for some $\beta\geq 1$, then one necessarily has 
\begin{equation}
  \label{eq:1}
d> q, \quad \beta\geq q^*,
\end{equation}
unless $b=0$. Indeed, this estimate is equivalent to the following:
\begin{equation}\label{eq:7}
\left\|D^2\,u\right\|_{\left(L^{\beta}(\Omega/\varepsilon)\right)^{d\times d}}+ \left\|\nabla u\right\|_{\left(L^{\beta}(\Omega/\varepsilon)\right)^d}\,\leq C_q\, \left\|f\right\|_{\left(L^{\beta}\cap L^q\right)(\Omega/\varepsilon)},
\end{equation}
for some constant $C$ independent of $\varepsilon\in (0,1)$ and of $\Omega$, where $u$ is a solution to
$-a_{ij}\partial_{ij} u + b_i\partial_i u = f$ in $\Omega/\varepsilon$ with, say, homogeneous Dirichlet boundary
condition. Next, let us consider the solution $u_\varepsilon$ of the problem
\begin{displaymath}
  a_{ij}\left(\frac x \varepsilon\right) \partial_{ij} u_\varepsilon + b_i \left(\frac x
    \varepsilon\right) \partial_i u_\varepsilon = f
\end{displaymath}
in $\Omega$, with homogeneous Dirichlet boundary conditions. In particular, this estimate implies
\begin{equation}\label{eq:12}
\left\|D^2\,u\right\|_{\left(L^{\beta}(\Omega/\varepsilon)\right)^{d\times d}}\,\leq C_q\, \left\|f\right\|_{\left(L^{\beta}\cap L^q\right)(\Omega/\varepsilon)}.
\end{equation}
Rescaling this equation, we find that
$v_\varepsilon(x) := u_\varepsilon\left(\varepsilon x\right)$ is solution to
\begin{displaymath}
  a_{ij}\partial_{ij}v_\varepsilon + \varepsilon b_i \partial_i v_\varepsilon = \varepsilon^2 f(\varepsilon x),
\end{displaymath}
in $\Omega/\varepsilon$. Hence, applying \eqref{eq:12}, we have 
\begin{multline*}
  \left\|D^2\,v_\varepsilon\right\|_{\left(L^{\beta}(\Omega/\varepsilon)\right)^{d\times d}}%+  \left\|\nabla
    % v_\varepsilon\right\|_{\left(L^{\beta}(\Omega/\varepsilon)\right)^d}
  \,\leq C_q\, \varepsilon^2
  \left\|f(\varepsilon\cdot)\right\|_{\left(L^{\beta}\cap L^q\right)(\Omega/\varepsilon)} \\ = C_q
  \left(\varepsilon^{2-\frac d q}\|f\|_{L^q(\Omega)} + \varepsilon^{2-\frac d {\beta}}\|f\|_{L^{\beta}(\Omega)}\right),
\end{multline*}
Going back to $u_\varepsilon$, this reads
\begin{equation}\label{eq:11}
  \varepsilon^{2-\frac d \beta} \left\|D^2 u_\varepsilon\right\|_{L^\beta(\Omega)} % + \varepsilon^{1-\frac d
    % \beta} \left\|\nabla u_\varepsilon\right\|_{L^\beta(\Omega)}
  \leq C_q
  \left(\varepsilon^{2-\frac d q}\|f\|_{L^q(\Omega)} + \varepsilon^{2-\frac d {\beta}}\|f\|_{L^{\beta}(\Omega)}\right),
\end{equation}
Finally, assuming that $a$ and $b$ are periodic, we apply standard homogenization technique to $u_\varepsilon$,
getting
\begin{displaymath}
  u_\varepsilon(x) = u^*(x) + \varepsilon \partial_j u^*(x) w_j\left(\frac x \varepsilon\right) + \varepsilon^2
  g\left(x,\frac x \varepsilon\right).
\end{displaymath}
Here, $w_j$ denotes the corrector associated to the above equation, and $u^*$ is the solution of the homogenized
problem, that is, the limit of $u_\varepsilon$ as $\varepsilon\to 0$. If all the data are smooth, we may assume that $g$ is smooth, hence,
\begin{displaymath}
  \partial_i u_\varepsilon(x) = \partial_i u^*(x) + \partial_j u^*(x) \partial_i w_j\left(\frac x \varepsilon\right) + O(\varepsilon),
\end{displaymath}
\begin{displaymath}
  \partial_{ik}u_\varepsilon(x) = \partial_{ik}u^*(x) + \frac 1 \varepsilon \partial_ju^*(x) \partial_{ik} w_j
  \left(\frac x \varepsilon\right) + O(1).
\end{displaymath}
These estimates imply that
\begin{equation}\label{eq:13}
  % \|\nabla u_\varepsilon\|_{L^\beta(\Omega)} = O(1), \quad
  \|D^2 u_\varepsilon\|_{L^\beta(\Omega)} \quad\text{scales as}\quad
  \frac 1 \varepsilon,
\end{equation}
unless $D^2 w_j = 0$, for all $j$ (recall that $u^*$ is independent of $w_j$). In the periodic case we are
studying here, this implies $\nabla w_j = 0$, that is, $b=0$. Inserting \eqref{eq:13} into \eqref{eq:11}, we find that
\begin{displaymath}
  \varepsilon^{1-\frac d \beta} \leq C  \left(\varepsilon^{2-\frac d q} + \varepsilon^{2-\frac d {\beta}}\right),
\end{displaymath}
where the constant $C$ depends on $a$, $b$, $f$, but not on $\varepsilon$. Letting $\varepsilon\to 0$, we thus find
$1-\frac d \beta \geq \min \left(2-\frac d q , 2-\frac d \beta\right)$, that is, \eqref{eq:1}.
\end{remark}

% \begin{remark}{\bf [XB: je laisse cette remarque car elle me semble juste, meme si c'est etrange vu la precedente.]}
% Another way of proving the statement of remark~\ref{rk:perte} is as follows:
% considering a solution $u$ of \eqref{eq:equation-adv-diff}, and defining $u_\lambda(x) = u(\lambda x)$, we have
% \begin{displaymath}
% - a_{ij}(\lambda x) \partial_{ij} u_\lambda + \lambda b_i(\lambda x) \partial_i u_\lambda = \lambda^2 f(\lambda x).  
% \end{displaymath}
% Hence, for $\lambda \leq 1$, we have a uniform bound on $\|a(\lambda \cdot)\|_{C^{0,\alpha}_{\rm unif}(\RR^d)}$ and $\|b(\lambda
% \cdot)\|_{C^{0,\alpha}_{\rm unif}(\RR^d)}$, hence the above estimate applies to $u_\lambda$. Thus,
% \begin{displaymath}
% \lambda^{2-\frac d \beta}\left\|D^2\,u\right\|_{\left(L^{\beta}(\RR^d)\right)^{d\times d}}+ \lambda^{1-\frac d
%   \beta}\left\|\nabla u\right\|_{L^{\beta}(\RR^d)}\,\leq C_q\, \lambda^{2-\frac d q}\left\|f\right\|_{\left(L^{q^*}\cap L^q\right)(\RR^d)}.
% \end{displaymath}
% Letting $\lambda \to 0$, we infer that $2-\frac d \beta\geq 2 -\frac d q$ and $1-\frac d \beta \geq 2 - \frac d
% q$. This implies \eqref{eq:1}.
% \end{remark}
\begin{remark}\label{rk:perte-2}
The case $b=0$ discussed in Remark~\ref{rk:perte} exactly corresponds to the condition
% In Remark~\ref{rk:perte}, the case $b=0$ is special and does not give \eqref{eq:1}. This is linked with the
% condition
given in \cite[Theorem B]{AL1991}, which states that for the equation
\begin{displaymath}
  -\div({\mathcal A}^{per} \nabla u) = f,
\end{displaymath}
an estimate of the form $\|D^2u\|_{\left(L^q(\RR^d)\right)^{d\times d}} \leq C \|f\|_{L^q(\RR^d)}$ can hold if and
only if $\div(\mathcal A^{per} ) = 0$. Actually, in the calculations of \cite[Section 3.1]{BLL-2017-1}, which are
recalled in Step 1 of the proof of Proposition~\ref{prop:non-divergence-form} above, we recover this
fact. The condition on $\mathcal A^{per}$ is equivalent to $\nabla w_j=0$ for all $j$, since
$-\div\left(\mathcal A^{per}(\nabla w_j + e_j)\right) = 0$, while $\div ({\mathcal A}^{per})  = \div (m_{per}
a^{per} - {\mathcal B}^{per})= m_{per} b^{per}$ shows that $\div\left(\mathcal A^{per}\right)=0$ if and only if $b^{per}=0$.
  % Tu pourrais aussi faire le lien avec la condition donnée par Avellaneda-Lin 1991, Theorem B, p 899. D'après les calculs rappelés dans notre premier article, Section 3.1, formules (36) et (43) (numérotation dans la version du 26 juillet), on a que l'estimée ne peut être bonne (dans le cas périodique) que si div ({\mathcal A}_{per})= 0, ce qui est exactement dire nabla w_p=0 pour tout p, puisque l'équation du correcteur se récrit : - div ({\mathcal A}_{per} (p+\nabla w_p))= 0. Or, div ({\mathcal A}_{per})  = div (m_{per} a^{per} - {\mathcal B}^{per})= m_{per} b^{per} , donc on retrouve la condition b=0 que tu donnes.
\end{remark}

\medskip

\begin{remark}
The norm of the right-hand side in \eqref{eq:estimee-informelle-non-div} is by definition
\begin{displaymath}
  \left\|f\right\|_{\left(L^{q^*}\cap L^q\right)(\RR^d)} = \left\|f\right\|_{L^q(\RR^d)} + \left\|f\right\|_{L^{q^*}(\RR^d)}.
\end{displaymath}
The presence of the second term is in fact necessary for the estimate \eqref{eq:estimee-informelle-non-div} to hold
true. Indeed, if $a_{ij}\partial_{ij}u=f$ and $D^2u\in L^{q^*}$
then $f\in L^{q^*}$. Moreover, we are going to use this estimate for functions belonging to spaces of the form
$L^q\cap L^\infty$. Having to consider functions that, in addition to being in $L^q$, belong to $L^{q^*}$, is not
a constraint in our setting. 
%{\bf [CLB: COMMENTER LE FAIT QU'IL Y A UNE INTERSECTION A DROITE DANS L'ESTIMEE : l'appartenance \`a $L^{q^*}$ est une condition "n\'ecessaire" (si $a_{ij}\partial_{ij}u=f$ et $D^2u\in L^{q^*}$ alors $f\in L^{q^*}$), et elle ne sera pas g\^enante dans la suite car le dual de $L^{q^*}\cap L^q$ est $L^{(q^*)'}+ L^{q'}$, et, comme  on manipulera seulement des fonctions qui sont aussi $L^\infty$, et que  $\left(L^{(q^*)'}+ L^{q'}\right)\cap L^\infty\subset L^{q'}\cap L^\infty$, on ne "verra" donc pas la partie $L^{(q^*)'}$.]
\end{remark}

% \medskip

% \noindent {\bf [CLB: Remarque perso: Tout cela \'etant, il y a quand m\^eme quelque chose de bizarre dans la mani\`ere dont j'\'enonce l'estim\'ee~\eqref{eq:estimee-informelle-non-div}: oui il y a une perte n\'ecessaire, mais c'est seulement si $b\not = 0$ (et encore ca doit d\'ependre des directions de $b$...) mais si $b=0$ elle est clairement sous-optimale, puisqu'il n'y a pas de perte en $1/d$ et qu'on n'est pas oblig\'e de supposer $q<d$ et $f\in L^q$ donc, pour obtenir $D^2 u\in L^{q^*}$ in fine ... il doit exister une autre formulation de l'\'enonc\'e du m\^eme r\'esultat qui soit optimale dans tous les cas... De m\^eme, pour l'\'enonc\'e du Corollaire concernant l'existence d'un correcteur, on \'enonce qu'il est sous-quadratique, mais dans le cas o\`u $b=0$ ou celui on peut mettre sous forme divergence, il est en fait sous-lin\'eaire...]} 

% \medskip

% \begin{remark}
% {\bf [COMMENTER LE CAS $d=2$ POUR LE COROLLAIRE~\ref{cor:corrector}: on a donc $s=1$, on a bien un correcteur, avec $D^2\tilde w_p\in L^\alpha$ pour $1/\alpha=1/r-1/2$. Au mieux $r=1$, d'o\`u $\alpha=2$. Pour autant, quoi ? ....]} 
% \end{remark}

% \medskip

\begin{remark}[On our assumptions $r$ and $s$ sufficiently small]
\label{rk:r-s-small}
Proposition~\ref{prop:non-divergence-form} and Corollary~\ref{cor:corrector} hold true under the assumption, in
particular, that  the perturbations $\tilde a$ and $\tilde b$ decay sufficiently fast to zero at infinity, namely
that they belong to $L^r$ and $L^s$ with $r$ and $s$ smaller than $d$.  Such a condition turns out to be,
qualitatively, necessary. And it is necessary not only to obtain a corrector with perturbation $\nabla\tilde w$ in
some $L^p$ space, but to obtain a corrector that satisfies the sharp condition to be imposed to a corrector, which
is only a \emph{consequence} of the condition $\nabla\tilde w\in L^p$: to be (strictly) sub-linear at infinity. In
order to show this is the case, we consider the simplistic one-dimensional situation where $a^{per}=1$, $\tilde
a=0$, $b=b^{per}+\tilde b$. The corrector equation then reads~$ -w''+b(1+w')=0$. The sufficient and necessary
condition for a periodic corrector $w_{per}$ to exist is~$\left\langle b^{per}\right\rangle=0$ (note that this  condition is indeed equivalent to
$\left\langle m_{per} b^{per}\right\rangle=0$ in this specific situation). That corrector is defined by~$w_{per}'= -1+\left\langle e^{B^{per}}\right\rangle^{-1} e^{B^{per}}$, where~$B^{per}=\int_0^xb^{per}(t)\,dt$. Then, any solution to the corrector equation reads (up to the addition of an irrelevant constant) as  $w=w_{per}+\tilde w$, where 
$$(\tilde w)'(x)=\left\langle e^{B^{per}}\right\rangle^{-1} e^{B^{per}(x)} (-1+e^{\tilde B(x)})$$ and~$({\tilde B})'=-\tilde b$.
 If we then impose on $\tilde w$ to be strictly sub-linear at $x=\pm\infty$, then we \emph{must} have~${\tilde
   B}(\pm\infty)=0$. In other words, both integrals~$\int_{-\infty}^{0}\tilde b(t)\,dt$
 and~$\int_{0}^{+\infty}\tilde b(t)\,dt$ must be well defined and~$\int_{-\infty}^{+\infty}\tilde
 b(t)\,dt=0$. %\textbf{[XB: attention, les integrales peuvent etre definies sans que la fonction soit $L^1$. Exemple
 %  $sin(x)/x$.]}
 %(condition qui est probablement un artefact de la dimension~1, puisque pour $d\geq 2$ et $\tilde b$ dans un $L^p$ convenable, on ne voit pas cette condition
It therefore in particular implies that $\tilde b$ has necessarily some integrability at infinity. For completeness, we check that the above conditions are indeed sufficient: the derivative~$(\tilde w)'$ then behaves as $\displaystyle \int_x^{+\infty}\tilde b(t)\,dt$ as $|x|\longrightarrow +\infty$, and~$\tilde w$ is  strictly sub-linear at infinity since ${{\tilde w(x)-\tilde w(0)}\over x}\approx \int_x^{+\infty}\tilde b(t)\,dt+\int_0^x {t\over x} \tilde b(t) dt.$ %{\bf [CLB: V\'erifier l'alg\`ebre des calculs]}.
%Il n'est pas s\^ur que cette solution sous-lin\'eaire soit telle $(\tilde w)'\in L^1$ (ce qui ne contredit pas la th\'eorie "g\'en\'erale" puisque $1$ n'est pas inf\'erieur \`a $d/2$ pour $d=1$). 
% On the other hand, if we only assume~$\tilde b\in L^p$ ($\cap L^\infty$), for some $p>1$, then  there is no reason for $\tilde w$ to be sub-linear at infinity. 
On the other hand, it is easy to build an example of $\tilde b\in L^p$ for some $p>1$, for which $\tilde b \notin
L^1$, and $\tilde w'$ grows exponentially at infinity. Think for instance of $\tilde b(x) \approx |x|^{-1/2}$ at
infinity, for which $\tilde w'(x) \approx e^{2 \sqrt{|x|}}$ at infinity.
\end{remark}

%\medskip
%
%\begin{remark}
%{\bf 
%\end{remark}
%
%\medskip

%%%%%%%%%%%%%%%%%%%
\section{Existence of the invariant measure}
\label{sec:adjoint}
%
%1) existence de la mesure invariante 
%
%2) avec en plus une remarque disant que si $b$ a une structure on peut faire diff\'eremment, c'est la Remark~\ref{rk:structure} (suivie de Remark~\ref{rk:structure-corrector} \`a la section suivante pour le correcteur)
%

%\medskip

We now consider the issue of existence (and uniqueness in a suitable class) of an invariant measure associated to equation~\eqref{eq:equation-adv-diff-epsilon}, that is  a positive function~$m$, actually bounded away from zero, $\inf m>0$, unique when appropriately normalized, solution to the  equation 
\begin{equation}
\label{eq:general-equation-adjoint}
-\partial_{i} (\partial_j (a_{ij}\,m) + b_i \,m)=0, 
\end{equation}
on $\RR^d$. We know from our previous study~\cite[Section~3]{BLL-2017-1} that this issue  is a straightforward consequence of the general estimate of the type~\eqref{eq:estimee-informelle-non-div}. The argument essentially goes by duality. 

\medskip

First of all, we know from the general theory (see e.g. \cite{blp}), that there exists a unique, periodic
measure~$m_{per}$, with normalized periodic average~$\left\langle m_{per}\right \rangle=1$, solution to 
\begin{equation}
\label{eq:general-equation-adjoint-periodic}
-\partial_{i} (\partial_j (a^{per}_{ij}\,m_{per} )+ b^{per}_i \,m_{per})=0. 
\end{equation}
Then we look for $m$ solution to \eqref{eq:general-equation-adjoint} as $m=m_{per}+\tilde m$ and rewrite \eqref{eq:general-equation-adjoint} as
\begin{equation}
\label{eq:general-equation-adjoint-tilde}
-\partial_{i} (\partial_j (a_{ij}\,\tilde m) + b_i \,\tilde m)=\partial_{i} (\partial_j (\tilde a_{ij}\,m_{per} )+ \tilde b_i \,m_{per}), 
\end{equation}
The key point for establishing the well-posedness of~\eqref{eq:general-equation-adjoint-tilde} 
is to show an {\it a priori} estimate on the solution to that equation. The conclusion follows by standard arguments made explicit in~\cite[Section~3]{BLL-2017-1}.
 
\medskip

Let us fix, as in Proposition~\ref{prop:non-divergence-form},  $1\leq r < d$, $\displaystyle 1\leq s< d$,
$1<q<d$. Recall our notation~$\displaystyle  {1\over {q^*}}={1\over q}- {1\over d}$. % and~$\displaystyle  {1\over {q^{**}}}={1\over q}- {2\over d}$
We denote by $q'$ the conjugate exponent of~$q$ and by $({q^*})'$ that of~${q^*}$. %, by $({q^{**}})'$ that
                                %of~${q^{**}}$. 
We have $\displaystyle {1\over q}+{1\over q'}=1$, $\displaystyle  {1\over {q^*}} + {1\over {(q^*)'}}=1$% , $\displaystyle  {1\over {q^{**}}} + {1\over {(q^{**})'}}=1$
. We consider the integral $\displaystyle \int\tilde m\,f$ for some arbitrary function $f\in \left(L^{q^*}\cap L^{q}\right)(\RR^d)$. 

Introducing the solution $u$ to $-a_{ij}\partial_{ij}u+b_j\,\partial_j u=f$ provided by Proposition~\ref{prop:non-divergence-form} and using~\eqref{eq:general-equation-adjoint-tilde}, we have 
\begin{eqnarray*}
\int\tilde m\,f&=&\int \tilde m\,\left(-a_{ij}\partial_{ij}u+b_j\,\partial_j u\right)\nonumber\\
&=&\int \left(-\partial_{i} (\partial_j (a_{ij}\,\tilde m) + b_i \,\tilde m)\right)\,u
\nonumber\\
&=&\int \left(\partial_{i} (\partial_j (\tilde a_{ij}\,m_{per} )+ \tilde b_i \,m_{per})\right)\,u\nonumber\\
&=&\int  m_{per}\,\left(\tilde a_{ij}\partial_{ij}u-\tilde b_j\,\partial_j u\right).
\end{eqnarray*}
The H\"older inequality and the estimate~\eqref{eq:estimee-informelle-non-div} successively yield
\begin{eqnarray}
\label{eq:adj2}
\left|\int\tilde m\,f\right|&\leq& 
\left\|m_{per}\right\|_{L^{\infty}(\RR^d)}
\left(
\left\|\tilde a\right\|_{L^{({q^*})'}(\RR^d)^{d\times d}}\,\left\|D^2 u\right\|_{L^{q^*}(\RR^d)^{d\times
                                   d}}\right. \nonumber\\
&& \left. +\left\|\tilde b\right\|_{L^{({q^{*}})'}(\RR^d)^d}\,\left\|\nabla u\right\|_{L^{q^{*}}(\RR^d)^d}\right)\nonumber\\
&\leq&C\,\left\|m_{per}\right\|_{L^{\infty}(\RR^d)}
\left(
\left\|\tilde a\right\|_{L^{({q^*})'}(\RR^d)^{d\times d}}+\left\|\tilde
       b\right\|_{L^{({q^{*}})'}(\RR^d)^d}\right)\nonumber\\
&&\times \left(\left\|D^2 u\right\|_{L^{q^*}(\RR^d)^{d\times d}} + \left\|\nabla u\right\|_{L^{q^{*}}(\RR^d)^d}\right)
\nonumber\\
&\leq&C\,\left\|m_{per}\right\|_{L^{\infty}(\RR^d)}
\left(
\left\|\tilde a\right\|_{L^{({q^*})'}(\RR^d)^{d\times d}}+\left\|\tilde
       b\right\|_{L^{({q^{**}})'}(\RR^d)^d}\right)\nonumber\\
&&\times \left\|f\right\|_{\left(L^{q^*}\cap L^{q}\right)(\RR^d)},
\end{eqnarray}
for some irrelevant constants $C$. 
By definition, 
%\begin{equation}
%\left\|\tilde m\right\|_{L^{q'}(\RR^d)}=\sup_{f\not=0\,\in L^{q}(\RR^d)}
%\frac{\left|\int\tilde m\,f\right|}{\left\|f\right\|_{L^{q}(\RR^d)}}
%\end{equation}
\begin{equation}
\label{eq:def-adj-norm}
\left\|\tilde m\right\|_{\left(L^{q'}+L^{(q^*)'}\right)(\RR^d)}=\sup_{f\not=0\,\in \left(L^{q}\cap L^{q^*}\right)(\RR^d)}
\frac{\left|\int\tilde m\,f\right|}{\left\|f\right\|_{\left(L^{q}\cap L^{q^*}\right)(\RR^d)}}
\end{equation}
We therefore infer from~\eqref{eq:adj2} and \eqref{eq:def-adj-norm} that $\tilde m\in \left(L^{q'}+L^{(q^*)'}\right)(\RR^d)$ with $1\leq q'<+\infty$, $1\leq (q^*)'<+\infty$,  provided $q$ is such that our assumptions~\eqref{eq:hyp1} on the integrability of $\tilde a$ and $\tilde b$ imply that 
$\tilde a\in \left(L^{({q^*})'}(\RR^d)\right)^{d\times d}$ and $\tilde b\in \left(L^{({q^{*}})'}(\RR^d)\right)^d$. This is the case when $r\leq ({q^*})'<+\infty$ and $s\leq ({q^{*}})'<+\infty$. The four conditions 
$$
\left\{
\begin{array}{l}
1\leq q'<+\infty,\\
1\leq (q^*)'<+\infty,\\
r\leq ({q^*})'<+\infty,\\
s\leq ({q^{*}})'<+\infty,
\end{array}
\right.
$$
reduce to 
%\begin{equation}
%\label{eq:q}
$\displaystyle {1\over q}\geq 1-\min\left({1\over r}-{1\over d}, {1\over s}-{1\over d}\right)$.
%\end{equation}
And we therefore obtain the best possible information on the integrability at infinity of $\tilde m$ when minimizing $q'$, that is maximizing $q$, that is taking an \emph{equality} in that equation, namely:
\begin{equation}
\label{eq:meilleur-q}
{1\over q}=1-\min\left({1\over r}-{1\over d}, {1\over s}-{1\over d}\right).
\end{equation}
On the other hand, we recall that by classical elliptic regularity results (see \cite[Theorem 9.11]{GT}),
$m_{per}\in C^{0,\alpha}(\RR^d)$ and $\tilde m \in C^{0,\alpha}(\RR^d)$. Moreover, standard results of periodic
homogenization \cite[Chapter 3, Section 3.3, Theorem 3.4]{blp} imply that $m_{per}$ is bounded away from $0$. Since $\tilde m\in L^q(\RR^d)$ and is
H\"older continuous, we have 
\begin{displaymath}
  \left\| \tilde m\right\|_{L^\infty(B_R^c)} \mathop{\longrightarrow}^{R\to +\infty} 0.
\end{displaymath}
Hence, for $R$ sufficiently large, 
\begin{equation}\label{eq:3}
  \forall x\in B_R^c, \quad m(x) \geq \frac12 \inf m_{per} >0.
\end{equation}
Applying the maximum principle on $B_R$, we deduce that $m\geq 0$ is valid in the whole space $\RR^d$. Next, we
apply Harnack inequality \cite{bogachev, bogachev-livre}, which implies that $m$ is bounded away from $0$. For the value of $q$ set in~\eqref{eq:meilleur-q}, we therefore obtain 
\begin{equation}
\label{eq:tilde-m-L-q}
\tilde m\in \left(L^{q'}\cap L^\infty\right)(\RR^d)\,\quad\hbox{\rm for}\quad {1\over q'}=\min\left({1\over r}-{1\over d}, {1\over s}-{1\over d}\right).
\end{equation}
% Finally, the measure $m$ is positive and bounded away from by a classical argument recalled in~\cite{BLL-2017-1}. {\bf [CLB: COMPLETER ICI. M\^eme commentaire que ci-dessus.]}
We collect our results in the following.
\medskip
\begin{corollary}
\label{cor:invariant-measure}
 %{\bf [CLB: On peut aussi placer l'\'enonc\'e avant la preuve, en changeant epsilon de la r\'edaction de ces deux pages.]} 
We assume \eqref{eq:aper+tildea}-\eqref{eq:hyp1} for some  $1\leq r<d$, $1\leq s<d$,
% donc $\max(r,s)<d$. 
and the condition~\eqref{eq:mper-bper}. Then there exists an invariant mesure $m$, solution to 
\eqref{eq:general-equation-adjoint}, that is
$-\partial_{i} (\partial_j (a_{ij}\,m) + b_i \,m)=0$. It reads as $m=m_{per}+\tilde m$, where $m_{per}$ is the
unique, normalized periodic invariant measure defined in~\eqref{eq:general-equation-adjoint-periodic}, and $\tilde
m$ belongs to $\left(L^{q'}\cap L^\infty\right)(\RR^d)$ where $q'$ is made precise in~\eqref{eq:tilde-m-L-q}. Such
a measure is unique, positive, bounded away from zero and H\"older continuous.
\end{corollary}

\medskip

% \begin{remark}
% \label{rk:2d-measure}
% {\bf [COMMENTER LE CAS $d=2$ POUR LE COROLLAIRE~\ref{cor:invariant-measure}: \eqref{eq:tilde-m-L-q} donne seulement , et ce "formellement" et en supposant $\tilde b\in L^1$ (ie $s=1$) que  $\tilde m\in L^\infty$ ...? ....]} Proposition d'argument rigoureux pour $d=2$. On prend $q=1$, d'o\`u $q^*=2$ (et $q^{**}=+\infty$). On a alors $D^2u\in L^2$. Si on a suppos\'e $\tilde b\in L\log L$ (au lieu de $\tilde b\in L^1$), alors par dualit\'e de $L\log L$ l'espace des fonctions expontentiellement int\'egrables $L_{exp}$  (voir [Adams, Theorem 8.27,
% p 277] et [Stein, Chapter 8] pour ces espaces et leurs propri\'et\'es) alors on obtient $\tilde m\in L^\infty$ (mais malheureusement pas dans $L^p$, $p<+\infty$ a priori. ...).
% \end{remark}

% \medskip

\begin{remark}[On coefficients with specific structure]
\label{rk:structure}
Our assumptions above are quite general. They apply without specific structure of the coefficients $a$ and $b$. If \emph{some} structure is assumed on these coefficients, then we suspect that the existence of an invariant measure may be proven using a different, more constructive approach. A simplistic example is $a^{per}_{ij}=\delta_{ij}$, $\tilde a_{ij}\equiv 0$ and $b$ (thus in particular $b^{per}$) is divergence-free.  Then we immediately observe that the periodic invariant measure is constant, and we normalize it to $m_{per}\equiv 1$, while $\tilde m\equiv 0$ (since we look for it in some $L^q(\RR^d)$). 
Similar examples may be constructed using different adequate coefficients $a$ and by "dividing"~$b$ by $a$. This
suffices to show that the presence of structure in the coefficients significantly changes the landscape. We wish to
concentrate here on an example which, although also simple, is more instructive. We again fix
$a^{per}_{ij}=\delta_{ij}$, $\tilde a_{ij}\equiv 0$, and this time set~$b^{per}_i\equiv 0$, and $\tilde
b=\nabla\tilde \psi$ for some function $\tilde \psi\in L^q(\RR^d)$ for some $1\leq q<+\infty$, $\psi$ sufficiently
regular (typically H\"older continuous, $C^{1,\alpha}$ so that the regularity assumed in~\eqref{eq:hyp1} is
satisfied). The perturbation $\tilde m$ of the periodic measure $m^{per}=1$ solves $\partial_j\left(\partial_j\tilde m+(1+\tilde m)\,\partial_j\tilde\psi\right)=0$. It is readily seen that $\tilde m=\exp{(-\tilde\psi)}-1$, so that the full invariant measure is $m=1+(\exp{(-\tilde\psi)}-1)=\exp{(-\tilde\psi)}$. Since $\tilde\psi$ vanishes at infinity (by regularity and integrability), $\tilde m$ behaves like~$\tilde\psi$ at infinity and also belongs to $L^q(\RR^d)$.  The point of this remark is that the exponent $1\leq q<+\infty$ may be arbitrarily large, in sharp contrast with both our "general" assumption $\tilde b\in \left(L^s(\RR^d)\right)^d$ for $s$ sufficiently small and our conclusion on $\tilde m\in L^\beta(\RR^d)$ again with $\beta$ small. Notice also that this observation does not contradict our considerations of Remark~\ref{rk:r-s-small}. Indeed, with this specific structure, $b^{per}\equiv 0$, $\tilde b=(\tilde\psi)'$ in our one-dimensional example there, and thus $(\tilde w)'=\exp{(-\tilde\psi)}-1$ does belong to $L^q(\RR^d)$.
\end{remark}
%%%%%%%%%%%%%%%%%%%
\section{Application to homogenization}
\label{sec:homog}

It is classical in the periodic case that the invariant measure allows one to recast (by multiplication) the
original problem as a problem for  an equation in divergence form. We have recalled the standard argument
in~\cite{BLL-2017-1} and above in Step 1 of the proof of Proposition~\ref{prop:non-divergence-form}. In the
present section, we extend it to the perturbed case with a drift and for simplicity we proceed in dimension $d\geq
3$. 

More precisely, we may rewrite~\eqref{eq:equation-adv-diff-epsilon}, and the associated corrector equation~\eqref{eq:corrector}, respectively as 
\begin{equation}
\label{eq:homog-div-1}
  -\div\left(\mathcal A_\varepsilon \nabla u^\varepsilon \right) = m_\varepsilon  f,
\end{equation}
and
\begin{equation}
\label{eq:corrector-div}
  -\div\left(\mathcal A (p+\nabla w_p) \right) = 0,
\end{equation}
with $ m_\varepsilon(x)= m(x/\varepsilon)$, with the elliptic matrix valued coefficient $\mathcal A_\varepsilon(x)=\mathcal A(x/\varepsilon)$ defined by  
\begin{equation}
\label{eq:homog-mathcal-1}
\mathcal A = m\, a - \mathcal B
\end{equation}
and the skew-symmetric matrix-valued coefficient~$\mathcal B$ is defined by
\begin{displaymath}
  \div(\mathcal B) = mb + \div(ma).
\end{displaymath}
Such a matrix may be proved to exist using the fact that $\div(mb + \div(ma)) = 0$, by definition of the measure
$m$. In the specific case of dimension $d=3$, we have
$$
%\begin{equation}
%\label{eq:mathcal-B}
  \mathcal B =\left(
  \begin{array}{ccc}
    0 & -B_3 & B_2 \\ B_3 & 0 & -B_1 \\ -B_2 & B_1 & 0
  \end{array}
  \right).
$$
%\end{equation}  
where the vector field~$B=\left(B_1,B_2,B_3\right)$ is defined by
%\begin{equation}
%\label{eq:B}
$\hbox{\rm curl}\, B=m\,b+\div(m\,a)$.
%\end{equation}
In our specific case, where $m=m_{per}+\tilde m$, $\mathcal B$ is defined as the sum $\mathcal B=\mathcal
B^{per}+\tilde{\mathcal B}$, where the periodic part $\mathcal B^{per}$ is obtained solving the periodic equation
$\div\,\mathcal B^{per}=m_{per}\,b^{per}+\div(m_{per}\,a^{per})$ (the right-hand side being divergence-free because of~\eqref{eq:general-equation-adjoint-periodic}, we recall) and where 
%\begin{equation}
%\label{eq:curlB}
\begin{equation}\label{eq:4}
\div\, \tilde{\mathcal B}=\tilde m\,b^{per}+\left(m_{per}+\tilde m\right)\,\tilde b+\div(\tilde m\,a^{per}+\left(m_{per}+\tilde m\right)\,\tilde a).
\end{equation}
% \end{equation}
The latter equation also has a divergence-free right-hand side by subtraction 
of~\eqref{eq:general-equation-adjoint-tilde} to \eqref{eq:general-equation-adjoint-periodic}. The matrix
$\tilde{\mathcal B}$, which is unique up to the addition of a constant, is found upon solving 
\begin{multline}
\label{eq:Delta-B}
-\Delta \tilde{\mathcal B}_{ij}=\partial_{jk}\Bigl(\tilde m a^{per}_{ik}+\left(m_{per}+\tilde m\right)\tilde
a_{ik}\Bigr) - \partial_{ik}\left(\tilde m a^{per}_{jk}+\left( m_{per}+\tilde m\right)\tilde a_{jk}\right) \\
 + \partial_j\left(\tilde m b_i^{per} + \left(m_{per}+\tilde m\right)\tilde b_i \right) - \partial_i\left(\tilde m b_j^{per} + \left(m_{per}+\tilde m\right)\tilde b_j  \right).
\end{multline}
Existence and uniqueness of the solution of this equation is proved using Calder\'on-Zygmund theory.
The detailed argument may be found in~\cite{BLL-2017-1} for the case $b=0$, with the result that $\tilde{\mathcal
  B}^{b=0} \in L^{q'}(\RR^d)$, with $\frac 1 {q'} = \min\left(\frac 1 r - \frac 1 d, \frac 1 s -
      \frac 1 d \right).$ In order to deal with $b$, since the
equation is linear, we only need to solve \eqref{eq:4} in the case $a=0$. For this purpose, we simply use the following
representation theorem:
\begin{multline*}
  \tilde {\mathcal B}_{ij}^{a=0} = (d-2)\frac {x_j}{|x|^d}*\left(\tilde m b_i^{per} + \left(m_{per}+\tilde
      m\right)\tilde b_i \right) \\- (d-2)\frac {x_i}{|x|^d}*\left(\tilde m b_j^{per} + \left(m_{per}+\tilde
      m\right)\tilde b_j \right).
\end{multline*}
Since $\frac {x_j}{|x|^d} \in L^{d/(d-1),\infty}(\RR^d)$ and $\tilde m b_i^{per} + \left(m_{per}+\tilde
      m\right)\tilde b_i \in L^{q'}(\RR^d)$, with $\frac 1 {q'} = \min\left(\frac 1 r - \frac 1 d, \frac 1 s -
      \frac 1 d \right),$ the Young-O'Neil inequality for Lorentz spaces \eqref{eq:young-oneil} implies that
    $\tilde{\mathcal B}^{a=0} \in L^\alpha(\RR^d),$ with ${1\over \alpha}= \frac{1}{(q')^*}=\min\left({1\over r}-{2\over d}, {1\over
        s}-{2\over d}\right).$ Finally, $\tilde{\mathcal B} = \tilde{\mathcal B}^{a=0} + \tilde{\mathcal B}^{b=0}$ satisfies
% Given the integrabilities \eqref{eq:hyp1} for $\tilde a \in \left(L^\infty(\RR^d)\cap L^r(\RR^d)\right)^{d\times d}$, $1\leq r<d$, and~$\tilde b \in \left(L^\infty(\RR^d)\cap L^s(\RR^d)\right)^d$, $1\leq s<d$, \eqref{eq:tilde-m-L-q} for $\displaystyle \tilde m\in \left(L^{q'}\cap L^\infty\right)(\RR^d)$, $\displaystyle  {1\over q'}=\min\left({1\over r}-{1\over d}, {1\over s}-{1\over d}\right)$, we observe that the limiting integrability (at infinity) in the right-hand side of~\eqref{eq:Delta-B} is that of the term $\tilde m\,b^{per}$ and therefore 
% we have  
\begin{equation}
\label{eq:B-tilde}
\tilde{\mathcal B}\in \left(L^\alpha(\RR^d)\right)^d\quad\hbox{\rm for}\quad {1\over \alpha}=\min\left({1\over r}-{2\over d}, {1\over s}-{2\over d}\right).
\end{equation}
We end up with the corrector problem \eqref{eq:corrector-div}, where
\begin{displaymath}
  \mathcal A = \underbrace{m_{per} a^{per} - {\mathcal B}^{per}}_{:= \mathcal A^{per}} + \underbrace{\tilde m
    a^{per} + \left(m_{per}+ \tilde m \right)\tilde a - \tilde {\mathcal B}.}_{:= \tilde{\mathcal A}} 
\end{displaymath}
The above considerations imply that $\tilde{\mathcal A}\in \left(L^\alpha\cap L^\infty\right)(\RR^d)$, with
$\alpha$ defined by \eqref{eq:B-tilde}. Hence, applying Proposition~2.1 of \cite{BLL-2017-1}, or Theorem~4.1 of
\cite{cpde-defauts}, we therefore find that the solution $w_p$ of \eqref{eq:corrector-div} exists, is unique up to
the addition of a constant, and reads $\nabla w_p = \nabla w_{p,per} + \nabla \tilde w_p$, where $w_{p,per}$ is the
periodic corrector associated with ${\mathcal A}^{per}$, and
\begin{equation}
  \label{eq:integrabilite_div}
  \nabla \tilde w_p \in L^\alpha(\RR^d)^d. 
\end{equation}
Compared to Corollary~\ref{cor:corrector}, we have seemingly lost some decay at infinity, since there, we have
$\nabla \tilde w_p \in L^{q^*}$, and $\frac 1 \alpha = \frac 1 {q^*} - \frac 1 d$, hence $\alpha >q^*.$

However, it is possible
to recover the fact that $\nabla \tilde w_p \in L^{q^*}(\RR^d)$ as follows: inserting $w_p = w_{p,per} + \tilde
w_p$ into \eqref{eq:corrector-div}, and using the fact that $-\div\left({\mathcal A}^{per}(\nabla
  w_{p,per}+p)\right)=0$, we write the equation satisfied by $\tilde w_p$:
\begin{displaymath}
  -\div\left((\mathcal{A}^{per}+ \tilde{\mathcal A})\nabla \tilde w_p \right) = \div\left(\tilde{\mathcal A} (\nabla w_{p,per}+p)\right).
\end{displaymath}
That is,
\begin{multline}\label{eq:2}
  -\div\left((\mathcal{A}^{per}+ \tilde{\mathcal A})\nabla \tilde w_p \right)
  =\div\left[(\tilde m a^{per} + \left(m_{per}+  \tilde m \right)\tilde a ) (\nabla w_{p,per}+p)\right]\\
 -\div\left[ \tilde{\mathcal B}\left( \nabla w_{p,per}+p\right)\right]
\end{multline}
Actually, the right-hand side of \eqref{eq:2} is exactly the right-hand side of \eqref{eq:corrector-bis} multiplied
by $m$. Hence, \eqref{eq:2} also reads 
\begin{equation}
  \label{eq:correcteur-modifie}
    -\div\left((\mathcal{A}^{per}+ \tilde{\mathcal A})\nabla \tilde w_p \right) = m \left( -\tilde b.p+\tilde a_{ij}\partial_{ij}w_{p,per}-\tilde b\cdot\nabla w_{p,per}\right).
\end{equation}
Next, we solve the following equation:
\begin{displaymath}
  -\Delta g = m \left( -\tilde b.p+\tilde a_{ij}\partial_{ij}w_{p,per}-\tilde b\cdot\nabla w_{p,per}\right) \in L^{\max(r,s)}(\RR^d),
\end{displaymath}
by defining
\begin{displaymath}
  g = \frac 1 {|x|^{d-2}} * m \left( -\tilde b.p+\tilde a_{ij}\partial_{ij}w_{p,per}-\tilde b\cdot\nabla w_{p,per}\right).
\end{displaymath}
Since $\nabla \frac 1 {|x|^{d-2}} \in L^{d/(d-1),\infty}(\RR^d)$, the Young-O'Neil inequality \eqref{eq:young-oneil} implies that
\begin{displaymath}
  \nabla g \in L^{q'}(\RR^d), \quad \frac 1 {q'} = \min\left(\frac 1 r - \frac 1 d, \frac 1 s - \frac 1 d \right).
\end{displaymath}
Hence, \eqref{eq:correcteur-modifie} also reads
\begin{displaymath}
    -\div\left((\mathcal{A}^{per}+ \tilde{\mathcal A})\nabla \tilde w_p \right) = \div\left(-\nabla g\right), \quad
      \nabla g\in L^{q'}(\RR^d).
\end{displaymath}
Applying Proposition~2.1 of \cite{BLL-2017-1}, we thus have $\nabla \tilde w_p\in L^{q'}(\RR^d)^d$. Thus, we
recover the result of Corollary~\ref{cor:corrector}.

Moreover, the fact that $\nabla \tilde w \in L^{q'}(\RR^d)$ allows to apply the theory of
\cite{josien,josien-these}, in order to find approximation results for the homogenization of
equation~\eqref{eq:equation-adv-diff-epsilon}.
We therefore find convergence theorems in $W^{1,p}$.

 %{\bf [CLB: v\'erifier et compl\'eter au besoin]} 

{\bf %[CLB: V\'erifier cet exposant $\alpha$  

}

% It follows from the above considerations that the matrix coefficient ${\mathcal A}$ defined in~\eqref{eq:homog-mathcal-1} and appearing in both equations~\eqref{eq:homog-div-1} and~\eqref{eq:corrector-div} is of the form ${\mathcal A}={\mathcal A}^{per}+\widetilde{\mathcal A}$ with ${\mathcal A}^{per}$ periodic defined in~\eqref{eq:mathcal-A-per} %=m_{per}\, a^{per} - {\mathcal B}^{per} $ 
% and 
% \begin{equation}
% \widetilde {\mathcal A}\,\in  \left(L^\alpha(\RR^d)\right)^{d\times d}\quad\hbox{\rm for the same}\quad \alpha\quad\hbox{\rm  as in}\quad \eqref{eq:B-tilde}
% \end{equation}
% In addition, because of the qualitative properties of $m$ and ${\mathcal B}$, ${\mathcal A}$ is elliptic. 

% \medskip

% We may therefore apply our results of~\cite{cpde-defauts}. They prove the existence of a corrector~$w_p$, solution to \eqref{eq:corrector-div} (or of course equivalently~\eqref{eq:corrector-div}). We find that ....{\bf [CLB: Compl\'eter et comparer avec le r\'esultat d'existence du correcteur obtenu directement
% ]}
% \medskip

\medskip

%\noindent {\bf [CLB: Quid de la fonction de Green, que peut-on prouver dessus?]}

Let us mention that, as pointed out in \cite{BLL-2017-1}, if $G$ is the Green function associated to
\eqref{eq:equation-adv-diff}, and if $\mathcal G$ is the Green function
associated to $-\div\left({\mathcal A}\nabla \cdot\right)$, we have
\begin{displaymath}
  \mathcal G (x,y) = m(y) G(x,y). 
\end{displaymath}
Therefore, all the estimates that are valid for the Green function $G$ yield adequate estimates on $\mathcal G$,
given the assumptions on $a,b$ and the regularity that they imply on $m$.
\medskip

\begin{remark}[Again on the case of coefficients with some specific structure]
\label{rk:structure-corrector}
We return here to the specific case we have examined in Remark~\ref{rk:structure}, that is  $a^{per}_{ij}=\delta_{ij}$, $\tilde a_{ij}\equiv 0$, $b^{per}_i\equiv 0$, and $\tilde b=\nabla\tilde \psi$ for some function $\tilde \psi\in L^q(\RR^d)$. We now look at the corrector functions. In this case, the corrector equation reads, for $p\in\RR^d$,  as~$-\Delta w_p+\nabla\tilde\psi\,.\,\nabla w_p=-p\,.\,\nabla w_p$. On the one hand, we evidently have $w_{p,per}=0$. On the other hand, multiplying the equation by the invariant measure $m=\exp{(-\tilde\psi)}$ yields $-\div\left(\exp{(-\tilde\psi)}\,(p+\nabla \tilde w_p)\right)=0$. Using our results on the equations in divergence form, we conclude to the existence of a corrector $\tilde w_p$ with $\nabla\tilde w_p\in L^q(\RR^d)$. Once again, we notice that $1\leq q<+\infty$ is arbitrary.
\end{remark}

%%%%%%%%%%%%%%%%%%%
% \section{Final remarks}
% \label{sec:final}

%%%%%%%%%%%%%%%%%%%%%%%%%%%%%%%%%%%%%
\section*{Acknowledgement} 
The work of the second author is partly
  supported by  ONR under Grant N00014-15-1-2777 and by EOARD, under Grant FA9550-17-1-0294. 
%{\bf [CLB: Mettre \`a jour l'autre article avec ce nouveau num\'ero de grant EOARD.]}
  
  \medskip

%{\bf [CLB: Ci-dessous dans la biblio, enlever \`a la fin les r\'ef\'erences non appel\'ees.]}
%%%%%%%%%%%%%%%%%%%%%%%%%%%%%%%%%%%%%%%%%%%%%%%%%%%%%%%%%%%%%%%%%%%%%%%%%%

\end{document}